\documentclass[12pt]{article}
\usepackage{amssymb}
\usepackage{amsmath}
\usepackage{amsthm}
\usepackage{color}
\usepackage{comment}
\usepackage{enumitem}		
\usepackage{geometry}		
\usepackage{graphicx}
\usepackage{authblk}
\usepackage{amssymb}
\usepackage{listings}
\usepackage{xcolor}
\usepackage{caption}
\usepackage{color}

\usepackage{subcaption}

\let\originalleft\left
\let\originalright\right
\renewcommand{\left}{\mathopen{}\mathclose\bgroup\originalleft}
\renewcommand{\right}{\aftergroup\egroup\originalright}

\newlist{romanlist}{enumerate}{3}
\setlist[romanlist]{label=\roman*),ref=(\roman*)}

\begin{document}

\newcommand{\cF}{\mathcal{F}}
\newcommand{\cR}{\mathcal{R}}
\newcommand{\cS}{\mathcal{S}}
\newcommand{\cT}{\mathcal{T}}
\newcommand{\cW}{\mathcal{W}}
\newcommand{\ee}{\varepsilon}
\newcommand{\rD}{{\rm D}}
\newcommand{\re}{{\rm e}}

\newcommand{\removableFootnote}[1]{\footnote{#1}}

\newtheorem{theorem}{Theorem}[section]
\newtheorem{corollary}[theorem]{Corollary}
\newtheorem{lemma}[theorem]{Lemma}
\newtheorem{proposition}[theorem]{Proposition}

\theoremstyle{definition}
\newtheorem{definition}{Definition}[section]
\newtheorem{example}[definition]{Example}

\theoremstyle{remark}
\newtheorem{remark}{Remark}[section]


\title{
On the analysis of a time varying noise-modulated heterogeneous coupled network of Chialvo neurons under the influence of electromagnetic flux.
}

\author[1]{Indranil~Ghosh}
\author[2, 1]{Sishu S.~Muni}
\author[3, 1]{Hammed O.~Fatoyinbo}

\affil[1]{School of Mathematical and Computational Sciences\\
Massey University\\
Colombo Road, Palmerston North, 4410\\
New Zealand\\
Email: i.ghosh@massey.ac.nz}
\affil[2]{Department of Physical Sciences\\
Indian Institute of Science and Educational Research Kolkata\\
Campus Road, Mohanpur, West Bengal, 741246\\
India\\
Email: sishu1729@iiserkol.ac.in}
\affil[3]{EpiCentre\\ School of Veterinary Science\\ Massey University, Palmerston North 4442, New Zealand\\
Email: h.fatoyinbo@massey.ac.nz}

\maketitle


\begin{abstract}
We perform a numerical study on the application of electromagnetic flux on a heterogeneous network of Chialvo neurons represented by a ring-star topology. Heterogeneities are realized by introducing additive noise modulations on both the central-peripheral and the peripheral-peripheral coupling links in the topology that not only vary in space but also in time. The variation in time is understood by two coupling probabilities, one for the central-peripheral connections and the other for the peripheral-peripheral connections respectively, that updates the network topology with each iteration in time. We have further reported the rich spatiotemporal patterns like two-cluster states, chimera states, traveling waves, coherent, and asynchronized states that arise throughout the network dynamics. We have also investigated the appearance of a special kind of asynchronization behavior called ``solitary nodes'' that have wide range of applications pertaining to real-world nervous systems. In order to characterize the behavior of the nodes under the influence of these heterogeneities, we have studied two different metrics called the ``cross-correlation coefficient'' and the ``synchronization error''. Additionally, to capture the statistical property of the network, for example, how complex the system behaves, we have also studied a measure called ``sample entropy". Various two-dimensional color-coded plots are presented in the study to exhibit how these metrics/measures behave with the variation of parameters. Finally, how the nodes synchronize or asynchronize is shown via one-dimensional bifurcation diagrams of the last instance of the main dynamical variable, i.e., the state variable associated with the membrane potential, against different network parameters.
\end{abstract}

\section{\label{sec:intro}Introduction}
Neurons form the fundamental units of the central and peripheral nervous system and supervise the mechanisms of complex information processing and responding to stimuli, by exchange of electrical and chemical signals between them. These complex dynamical behaviors exhibited by neurons can be represented and studied with the help of dynamical systems \cite{izhikevich2007dynamical, ibarz2011map} like ordinary differential equations and maps, leading to the science of neurodynamics. Recently, neurodynamics has become an emerging field of research and has attracted a lot of attentions \cite{izhikevich2007dynamical} from mathematicians, biologists, computer scientists, to name a few. These dynamical systems-oriented models, which mimic many neuronal behaviors, have been confirmed experimentally \cite{HGu13} too. Examples of models that are represented as continuous dynamical systems include Hodgkin-Huxley  model \cite{HoHu52}, Hindmarsh-Rose system \cite{hindmarsh1984model}, FitzHugh-Nagumo neuron system \cite{fitzhugh1961impulses}, etc.,  whereas examples for the ones represented as discrete systems involve Rulkov neuron system \cite{rulkov2002modeling} and Chialvo neuron system \cite{chialvo1995generic}. Scarcely any research attention is given to the study of neurodynamical models that are portrayed as discrete maps. Motivated, we focus on an improved model of a network of Chialvo neurons with heterogeneities incorporated, which we believe is a good imitation of a real-world nervous system. It is important to study the corresponding dynamical behaviors to gain an insight about how a nervous system might behave in reality.\par
One of the striking features exhibited by an ensemble of neurons is the phenomenon of synchronization. Synchronization is a universal concept in dynamical systems that have been studied in fields ranging from biology, to physics, to neuroscience, to economics \cite{luo2009theory, sumpter2006principles, boccaletti2018synchronization, pikovsky2002synchronization, strogatz2012sync}. Relevant to our study in neurodynamics, a synchronized state can both pertain to a normal or abnormal cognitive state \cite{uhlhaas2006neural}. In the latter case, it is of utmost importance to study and understand the dynamics of neurons that are not completely synchronized, as has been mentioned in the paper \cite{rybalova2019solitary}. Two of the interesting states that represent asynchronicity in the neural nodes are \textit{chimera states} and \textit{solitary nodes}. A ``chimera'' state \cite{kuramoto2002coexistence, panaggio2015chimera, abrams2004chimera, strogatz2000kuramoto, jaros2015chimera}, discovered by Kuramoto, is characterized by the simultaneous existence of coherent and incoherent nodes from a specific ensemble after a precise time. Chimera states have been reported in various natural phenomenon like sleeping patterns of animals \cite{rattenborg2000behavioral, mathews2006asynchronous, glaze2021neural}, flashing of fireflies \cite{haugland2015self}, and many more \cite{wickramasinghe2013spatially, panaggio2015chimera, martens2013chimera}. On the other hand, nodes falling under the ``solitary" regime \cite{rybalovasolitary, rybalova2021interplay, semenova2018mechanism, maistrenko2014solitary, bukh2017new, rybalova2019solitary, shepelev2020role} are the ones which behave completely different from the coherent nodes in a particular ensemble at a specific time. They get separated from the main typical cluster and oscillate with their own particular amplitude. Additionally, traveling waves are reported in various neuron systems \cite{baker2022traveling, ermentrout2001traveling}. The authors have reported traveling waves in heterogeneous neuron systems and have claimed that it is responsible for many features of cortical dynamics, spiking time variability, and fluctuations in membrane potential. There have been research works concerning networks of dynamical systems in which the nodal behaviors were found to be rich. It is interesting to note that in a network of neurodynamical systems, the researchers have found chimera states, cluster synchronization, traveling waves, solitary states, etc. A variety of spatiotemporal patterns are reported in the heterogeneous network considered in this study,  such as solitary states, traveling waves, synchronized state, asynchronized states, etc.\par
The fact that the dynamical behaviors portrayed by neurons are remarkably complex, demands the requirement of statistical tools to study and quantify their complexity. Measures like \textit{spatial average density}, and \textit{entropy} are effective tools to study the complexity in the dynamics of neurons. Banerjee and Petrovskii \cite{BanPet2011} applied the spatial average density approach to show that the spatial irregularity in an ecological model is indeed chaotic. Similar approach is considered in the reference \cite{BanBan2012}. Entropy is an ubiquitous concept first introduced by C. E. Shannon in his revolutionary work \cite{shannon2001mathematical}. Since then, it has had a widespread application in various domains of research, including neuroscience \cite{zbili2021quick, timme2018tutorial}. In another study \cite{wang2019sample} authors report the study of applying entropy to EEG data from patients and report that Alzheimer’s disease could result in complexity loss. More relevant paper  \cite{gomez2010entropy} related to this topic can also be found. Other applications of entropy in neuroscience include study of topological connectivity in neuronal networks \cite{vicente2011transfer, ursino2020transfer, ito2011extending}, and estimation of the upper limit on the information content in the action potential of neurons \cite{street2020upper}. Motivated by these, we utilize the tool of \textit{sample entropy} \cite{richman2000physiological} to analyze time-series data of the spatially averaged action potential generated from our model and study the complexity of the dynamics.\par
Empirical evidence says that the functioning of neurons are affected by many external factors such as temperature \cite{buzatu2009temperature}, light \cite{orlowska2021neuronal}, electromagnetic radiations \cite{muni2022discrete, discChialvo, Ham_et_al22} and many more. For example, signature of growth in embryo neuronal cells was observed with the application of electromagnetic flux \cite{kaplan2016electromagnetic}. These evidences motivate us to perform a study on the effect of external electromagnetic flux in a network of Chialvo neurons. The ring-star network topology \cite{muni2020chimera} serves as a promising candidate. Most recently, the dynamical effects of discrete Chialvo neuron map under the influence of electromagnetic flux have been studied \cite{discChialvo}, where after understanding a single neuron system, analysis on a ring-star network of the neurons has been conducted, making it flexible to work with all the possible scenarios of ring, star, and ring-star topologies. The network that has been mathematically set up in the work, consists of $N$ Chialvo neurons arranged in the formation of a ring-star topology, where the central node is connected to all the $N-1$ nodes with the homogeneous coupling strength $\mu$ and the peripheral nodes are connected to each of their $R$ neighbors with a coupling strength $\sigma$. The nodes follow a periodic boundary condition, meaning that the $(N-1)^{th}$ node is connected to the $2^{nd}$ node to complete the loop. Although the studies in the paper have exhibited some rich dynamics, the topology of the network is too simple to be close to a real-world connection of the nervous system. The study of these neurodynamical systems on this simple network topology gave further motivation to consider even more complex topologies to mimic the real-world connection between the neurons in the brain. Extending on the work, we have tried to tackle this issue by studying the spatiotemporal dynamics on introduction of heterogeneous and time-varying coupling strengths to the network topology of the ring-star Chialvo neurons under the influence of an electromagnetic flux.\par
In this paper, we mainly catalogue the presence of several interesting temporal phenomena i.e., solitary nodes, chimera, and traveling waves, to name a few, exhibited by the nodes in a heterogeneous network of Chialvo neurons. Heterogeneity is introduced both in space and time. In another study \cite{buscarino2015chimera}, the authors have reported chimera states in time-varying links of a network formed from two coupled populations of Kuramoto oscillators. To quantify what solitary states mean in our system model, we work with a metric called ``cross-correlation'' coefficient \cite{vadivasova2016correlation}. This metric lets us decide the regime a node lies in.\par
The goals of this paper are as follows:
\begin{enumerate}
    \item  Improve the ring-star topology of the Chialvo neuron network introduced in the original paper \cite{discChialvo}, by inclusion of heterogeneous links between nodes that vary not only in space but also in time, depending on noise modulations and specific coupling probabilities,
    \item Report the appearance of rich spatiotemporal patterns throughout the temporal dynamics of the heterogeneous network,
    \item Study the emergence of an important asynchronous behavior called \textit{solitary nodes} and characterize it using quantitative measures like \textit{cross-correlation coefficient} and \textit{synchronization error} that we set up according to our network topology,
    \item Use  the statistical tool called \textit{sample entropy} on the time series data of spatially averaged action potential, to gain an intuition on the extent of complexity, and
    \item Developing an understanding of the innate mechanism of whether and how the various measures relate by employing numerical simulations.
\end{enumerate}

We organize the paper in the following way: in Sec.~\ref{sec:model}, we introduce our improved heterogeneous model for a ring-star network of Chialvo neurons, followed by setting up of the quantitative metrics like cross-correlation coefficient and synchronization error in Sec.~\ref{sec:metrics}. Next, in Sec.~\ref{sec:stat_ts} we establish how we are going to apply the measures of sample entropy and maximum Lyapunov exponents on time series data that get generated from our model dynamics. Then in Sec.~\ref{sec:results} we show results from running various numerical experiments (plotting phase portraits, spatiotemporal patterns, recurrence plots, time series plots, several two regime color plots, bifurcation plots for synchronization, etc.), perform time series analysis, and try drawing inferences about the behavior and dynamics of our heterogeneous neuronal network model. Note that all simulations are performed in \texttt{Python}. Finally, in Sec.~\ref{sec:conclusions} we provide concluding remarks and future directions.

\section{\label{sec:model}System Modelling}
The two-dimensional discrete map proposed by Chialvo \cite{chialvo1995generic} in 1995 that corresponds to the single neuron dynamics is given by
\begin{align}
\label{eq:chialvoo_original}
        x(n+1) &= x(n)^{2} e^{(y(n) - x(n))} + k_{0},\\
    y(n+1) &= a y(n) - b x(n) + c,
\end{align}
where $x$ and $y$ are the dynamical variables representing activation and recovery respectively. Two of the control parameters are $a<1$ and $b<1$ which indicate the time constant of recovery and the activation dependence of recovery respectively during the neuron dynamics. Control parameters $c$ and $k_0$ are the offset and the time independent additive perturbation. Throughout this study we have kept $a = 0.89, b = 0.6, c = 0.28,$ and $k_0 = 0.04$.

The model was improved \cite{discChialvo} to a three-dimensional smooth map by inclusion of electromagnetic flux, realized by a memristor. The improved system of equations is give by

\begin{align}
\label{eq:chialvo_improved}
x(n+1)&= x(n)^{2} e^{(y(n)-x(n))} + k_{0} + k x(n) M(\phi(n)),\\
y(n+1) &= a y(n) - bx(n) + c,\\
\phi(n+1) &= k_{1}x(n) - k_{2}\phi(n),
\end{align}
where $M(\phi) = \alpha + 3 \beta \phi^2$ is the commonly used memconductance function \cite{chua2012hodgkin} in the literature with $\alpha, \beta, k_1$ and $k_2$ as the electromagnetic flux parameters. In \eqref{eq:chialvo_improved}, the $k_{1}x$ term denote the membrane potential induced changes on magnetic flux and the $k_{2}\phi$ term denote the leakage of magnetic flux \cite{discChialvo}. The parameter $k$ denotes the electromagnetic flux coupling strength with $kxM(\phi)$ as the induction current. We note that $k$ can take both positive and negative values. Throughout the paper we have used $\alpha= 0.1, \beta = 0.2, k_1 = 0.1, k_2 = 0.2$. We also restrict $k$ in the range $[-1, 4]$. \par

Recently \cite{discChialvo}, a ring-star network of the Chialvo neurons under the effect of electromagnetic flux has been considered where the coupling strength is set to be homogeneous, and they do not vary with time. We further improve the ring-star model with the incorporation of heterogeneous coupling strengths $\mu_m(n)$ and $\sigma_m(n)$. The improved version is then given by 

\begin{align}
    \label{eq:ChialvoDiscrete}
    x_m(n+1)&= x_m(n)^2e^{y_m(n)-x_m(n)} + k_0 + kx_m(n)M(\phi_m(n)) \nonumber \\
    &+\mu_m(n)(x_m(n)-x_1(n)) \nonumber \\
    &+\frac{1}{2R}\sum_{i=m-R}^{m+R}\sigma_i(n)(x_i(n)-x_m(n)), \\ 
    y_m(n+1)&=ay_m(n)-bx_m(n)+c, \\ 
    \phi_m(n+1) &= k_1x_{m}(n)-k_2\phi_m(n),
\end{align}
whose central node is further defined as
\begin{align}
    \label{eq:ChialvoDiscreteCentral}
    x_1(n+1) &= x_1(n)^2e^{(y_1(n)-x_1(n))}+k_0+kx_1(n)M(\phi_1(n)) \nonumber \\
            &+ \sum_{i=1}^N\mu_i(n)(x_i(n)-x_1(n)),  \\
    y_1(n+1) &= ay_1(n) - bx_1(n)+c, \\
    \phi_1(n+1) &= k_1x_1(n) - k_2\phi_1(n), 
\end{align}
having the following boundary conditions:

\begin{align}
    \label{eq:BoundaryConditions}
    x_{m+N}(n)= x_m(n), \\ 
    y_{m+N}(n) = y_m(n), \\ 
    \phi_{m+N}(n) = \phi_m(n),
\end{align}
where $n$ represents the $n^{\rm th}$ iteration, and $m = 1, \ldots, N$ with $N$ as the total number of neuron nodes in the system. For the sake of making the model more complex and closer to a realistic nervous system we introduce heterogeneities to the coupling strengths $\sigma_m(n)$ and $\mu_m(n)$ both in space and time. In space, the heterogeneities are realised following the application of a noise source with a uniform distribution \cite{rybalova2022response, vadivasova2016correlation} given by
\begin{align}
\label{hetero_coup}
\sigma_m(n) &= \sigma_0 + D_{\sigma}\xi_{\sigma}^{m, n}, \\ 
\mu_m(n) &= \mu_0 + D_{\mu}\xi_{\mu}^{m, n},
\end{align}
where $\sigma_0$ and $\mu_0$ are the mean values of the coupling strengths $\mu_m$ and $\sigma_m$ respectively. Throughout the paper we keep $\sigma_0 \in [-0.01, 0.01]$ and $\mu_0\in [-0.001, 0.001]$. The noise sources $\xi_{\sigma}$ and $\xi_{\mu}$ for the corresponding coupling strengths are real numbers randomly sampled from the uniform distribution $[-0.001, 0.001]$. Finally, the $D$'s refer to the ``noise intensity" which we restrict in the range $[0, 0.1]$. Note that we have used negative (inhibitory) coupling strengths in this study. They represent a significant proportion of neuronal connectivity in the human nervous system. The authors in the papers \cite{tsigkri2018synchronization, tsigkri2016multi} have mentioned them in their course of simulations of the leaky Integrate-and-Fire (LIF) model. Also, negative coupling strengths are included via a rotational coupling matrix (See Eq. (2) in the paper \cite{omelchenko2013nonlocal}) during simulations of FitzHugh-Nagumo neuronal models. Heterogeneity in time is introduced by considering the network having time-varying links \cite{kohar2014synchronization, de2015effect, mondal2008rapidly} depending on the two coupling probabilities $P_{\mu}$ and $P_{\sigma}$, which govern the update of the coupling topology with each iteration $n$. The probability with which the central node is connected to all the peripheral nodes at a particular $n$ is denoted by $P_{\mu}$. Likewise, the probability with which the peripheral nodes are connected to their $R$ neighboring nodes is given by $P_{\sigma}$. We have studied the spatiotemporal dynamics of the improved topological model under the variation of the seven most important parameters which are $k, \sigma_0, \mu_0, D_{\sigma}, D_{\mu}, P_{\sigma}$ and $P_{\mu}$.
\begin{figure}
    \centering
    \includegraphics[scale=1.5]{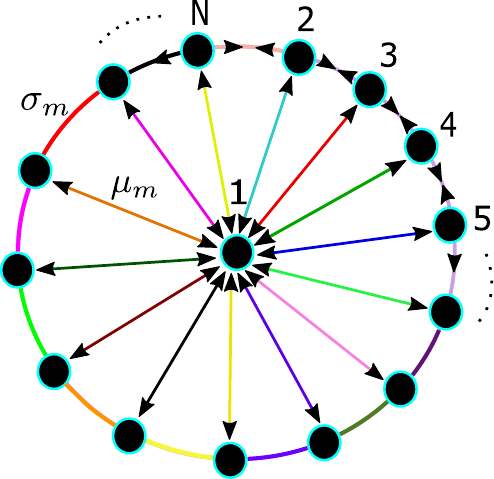}
     \caption{Heterogeneous ring-star network of Chialvo neurons. The nodes are numbered from $1, \ldots, N$. The star and ring coupling strengths are denoted by $\mu_{m}$ and $\sigma_{m}$ for each node $m=1, \ldots, N$ respectively. Different colors in the ring-star topology signify a range of heterogeneous values of $\mu_m$ and $\sigma_m$.}
 \label{fig:StochasticRingStar}
\end{figure}

\section{\label{sec:metrics}Quantitative metrics and time series analysis}
To quantize the solitary nodes and determine the extent of synchronization in the system after a satisfactory number of iterations, we employ two metrics known as the \textit{cross-correlation coefficient} \cite{semenova2018mechanism, vadivasova2016correlation} and the \textit{synchronization error} \cite{mehrabbeik2021synchronization}. We also realise the complexity of the time series data that get generated from our simulations using a measure known as s\textit{sample entropy}.

\subsection{\label{sec:cross_coeff}Cross-correlation coefficient}
Keeping in mind that our system has a ring-star topology, we must define the coefficient in such a way that it captures the correct collective behaviour of the network dynamics. The general definition of the cross-coefficient denoted by $\Gamma_{\rm i,m}$ is given by
\begin{align}
\label{ccc}
    \Gamma_{\rm i,m} = \frac{\langle \tilde{x}_i(n) \tilde{x}_m(n)\rangle}{\sqrt{\langle(\tilde{x}_i(n))^2 \rangle \langle(\tilde{x}_m(n))^2 \rangle}}.
\end{align}
The \textit{averaged cross-correlation coefficient} over all the units of the network is given by,
\begin{align}
\label{gcc}
    \Gamma = \frac{1}{N-1}\sum_{m=1, m \neq i}^N \Gamma_{\rm i,m}.
\end{align}
Throughout the paper, we use $\Gamma_{2,m}$, denoting the degree of correlation between the first peripheral node of the ring-star network and all the other nodes, including the central node. The average is calculated over time with transient dynamics removed and $\tilde{x}(n) = x(n) - \langle x(n) \rangle$ refers to the variation from the mean. Note that $\langle \rangle$ denotes the average over time. Everywhere in the paper we take $20000$ iterations in time, of which we discard the first $10000$ points for the dynamical variables, such that the transient dynamics is discarded, and thereafter perform all the calculations and simulations. We use $\Gamma_{\rm 2,m}$ to characterize all the necessary regimes for our study. When the nodes have $\Gamma_{\rm 2,m} \sim 1$, they lie in the \textit{coherent domain}. Due to the noise modulations, there might be node behaviors that are realized by their $\Gamma_{\rm 2,m}$ values lying in the range $[-0.15, 0.75]$. Finally, the solitary regime is characterized by the domain where the nodes have $\Gamma_{\rm 2,m} \in [-0.38, -0.15]$. Note that these values were selected after running numerous simulations confirming the fact that it is the local dynamics of the network of oscillators that govern the characteristics of solitary nodes from system to system. $\Gamma \to 1$ implies that almost all the nodes are clustered in the coherent regime after a specific time without any transient dynamics.  

\subsection{\label{sec:sync_error}Synchronization error}
The averaged \textit{synchronization-error} for the nodes in a system is given by
\begin{align}
\label{eq:sync_error}
    E = \frac{1}{N-1} \sum_{m=1, m \neq 2}^{N} \langle | x_2(n) - x_m(n) | \rangle,
\end{align}
where we again consider node number $N=2$ as the baseline for the calculation as has been done in the case of \textit{cross-correlation} coefficient, and $n$ represents the $n^{\rm th}$ iteration. Note that $E \to 1$ implies that the nodes in the system are moving towards an asynchronous behavior, with $E=1$ depicting a complete asynchrony. Similarly $E \to 0$ implies synchronization. Like the first metric, $\langle \rangle$ denotes the average over time.

\subsection{\label{sec:stat_ts}Sample entropy: a measure for complexity}
Additionally, we perform a statistical analysis of the dynamics of our system to determine how complex it is. In order to do so, we take the spatial average of all the $N$ nodes at a particular time $n$ and generate a time series data for the action potential, which we utilize to calculate the \textit{sample entropy (SE)}. SE tells us about how complex the time series data is. A high value indicates unpredictability in the behavior, thus offering more complexity \cite{he2021discrete, hansen2017sample}. To calculate  SE from the time series data, we utilize an open-source package called \texttt{nolds} \cite{nolds}, that provides us with the function \texttt{nolds.sampen()}. 

\section{Results}
\label{sec:results}

Every time we run a computer experiment, we use a pseudo-random generator to initialize the action potential $x$ between $0$ and $1$. Additionally $y$ and $\phi$ are initialized to $1$. This is done for all the $N=100$ nodes. In each simulation, the total number of iterations is set to $20000$, of which we discard the first $10000$ to remove transient dynamics. First, we briefly cover the single neuron dynamics under the fixed parameter values as mentioned in the previous sections, before moving to the analysis of the network.

\subsection{Single neuron dynamics}
\begin{figure}
    \centering
    \includegraphics[scale=.5]{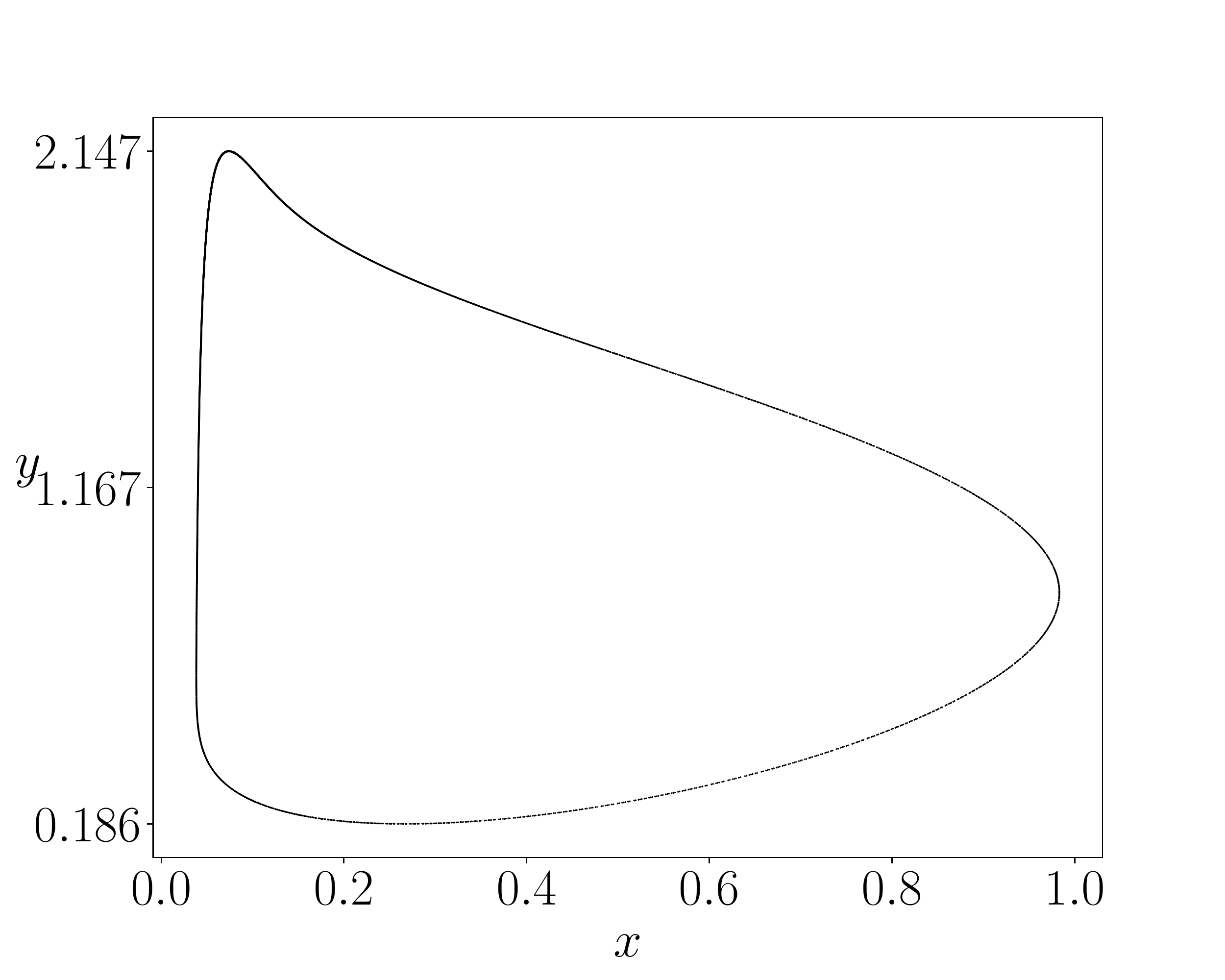}
     \caption{Phase portrait of a single neuron dynamics. Parameters are set to be $a = 0.89, b = 0.6, c = 0.28, k_0 = 0.04, \alpha= 0.1, \beta = 0.2, k_1 = 0.1,$ and $k_2 = 0.2$. A closed invariant curve (quasiperiodicity) is observed.}
 \label{fig:single pp}
\end{figure}

\begin{figure}
    \centering
    \includegraphics[scale=.08]{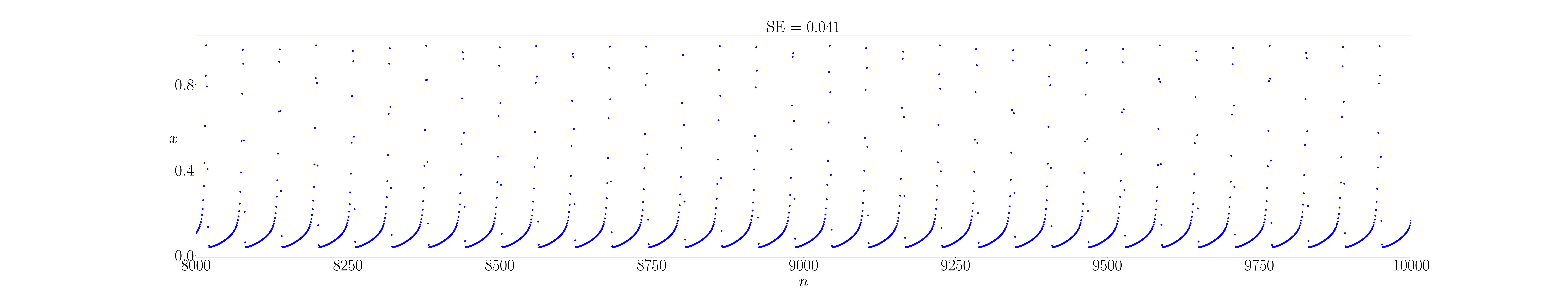}
     \caption{Time series plot for the corresponding neuron. The sample entropy is $0.041$, which is quite low. This indicates less disorderedness in the dynamics as can also be seen from not so irregular behaviour in the time series.}
 \label{fig:single TS}
\end{figure}

For a single neuron, we set the parameter values to $a = 0.89, b = 0.6, c = 0.28, k_0 = 0.04, \alpha= 0.1, \beta = 0.2, k_1 = 0.1,$ and $k_2 = 0.2$. Additionally, we set $k=-1$. The corresponding phase portrait and the time series are given in Fig.~\ref{fig:single pp} and Fig.~\ref{fig:single TS} respectively. The sample entropy value is calculated to be $0.041$, indicating quite a low complexity. Looking at the phase portrait, we observe that the attractor is a closed invariant curve. Using the $QR-$factorization method as was done recently \cite{discChialvo}, it can be noted that the maximum Lyapunov exponent is $\sim 0$, exhibiting a quasi-periodic dynamics for the fixed set of parameter values. Keeping this in mind, we  analyze the system of $N=100$ such neurons arranged in ring-star topology in the next sections. Note that for the above selected local dynamical parameters $a, b, c, k_0, \alpha= , \beta, k_1,$ and $k_2$, we perform the whole analysis throughout. In case, these parameters are set different, the dynamical behaviors of the single neuron and the network of neurons are expected to change accordingly.

\subsection{Phase portraits, spatiotemporal patterns and recurrence plots}

\begin{figure*}
    \centering
    \includegraphics[scale=0.7]{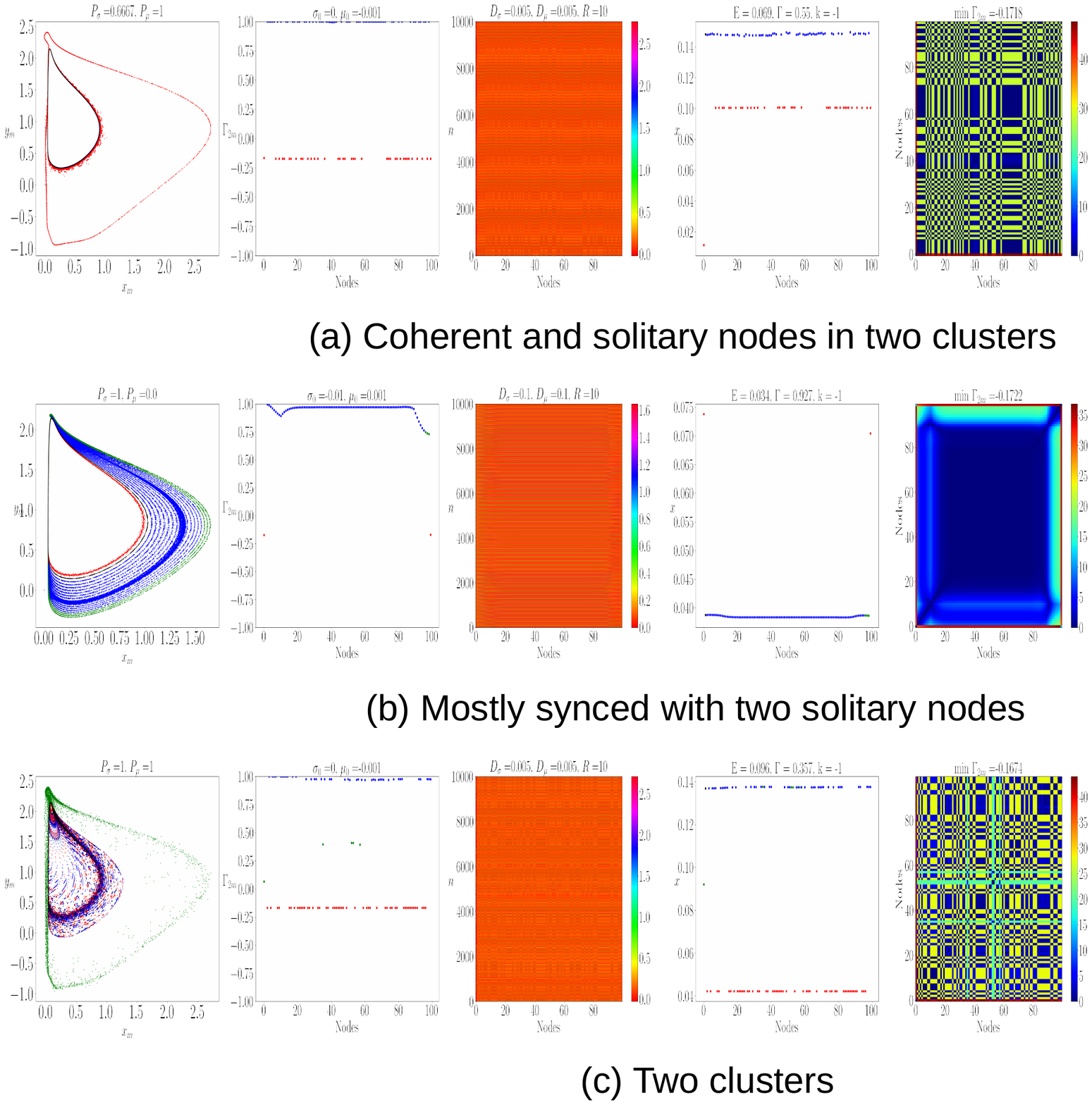}
     \caption{From left to right: phase portrait, $\Gamma_{\rm 2, m}$ plot for each node, spatiotemporal diagram, last instance plot of $x_m$ for each node, and recurrence plot. Parameters are set as (a) $\sigma_0 = 0, \mu_0=-0.001, D_{\sigma}=0.005, D_{\mu}=0.005, P_{\sigma}=0.66666, P_{\mu}=1,$ and $k=-1$, (b) $\sigma_0 = -0.01, \mu_0=0.001, D_{\sigma}=0.1, D_{\mu}=0.1, P_{\sigma}=1, P_{\mu}=0,$ and $k=-1$, and (c) $\sigma_0 = 0, \mu_0=-0.001, D_{\sigma}=0.005, D_{\mu}=0.005, P_{\sigma}=1, P_{\mu}=1,$ and $k=-1$. The normalized cross-correlation coefficient and the normalized syncronization error values are (a) $\Gamma = 0.55,$ $E = 0.069$, (b) $\Gamma = 0.927,$ $E = 0.034$, and (c) $\Gamma = 0.357,$ $E = 0.096$. In the first, second and fourth plots, the nodes lying in the solitary regime are denoted by red dots, the nodes with $\Gamma_{\rm 2,m} \in [-0.15, 0.75]$ are denoted by green dots and the nodes in the coherent domain are denoted by blue dots.}
 \label{fig:5plots1}
\end{figure*}

\begin{figure*}
    \centering
    \includegraphics[scale=0.7]{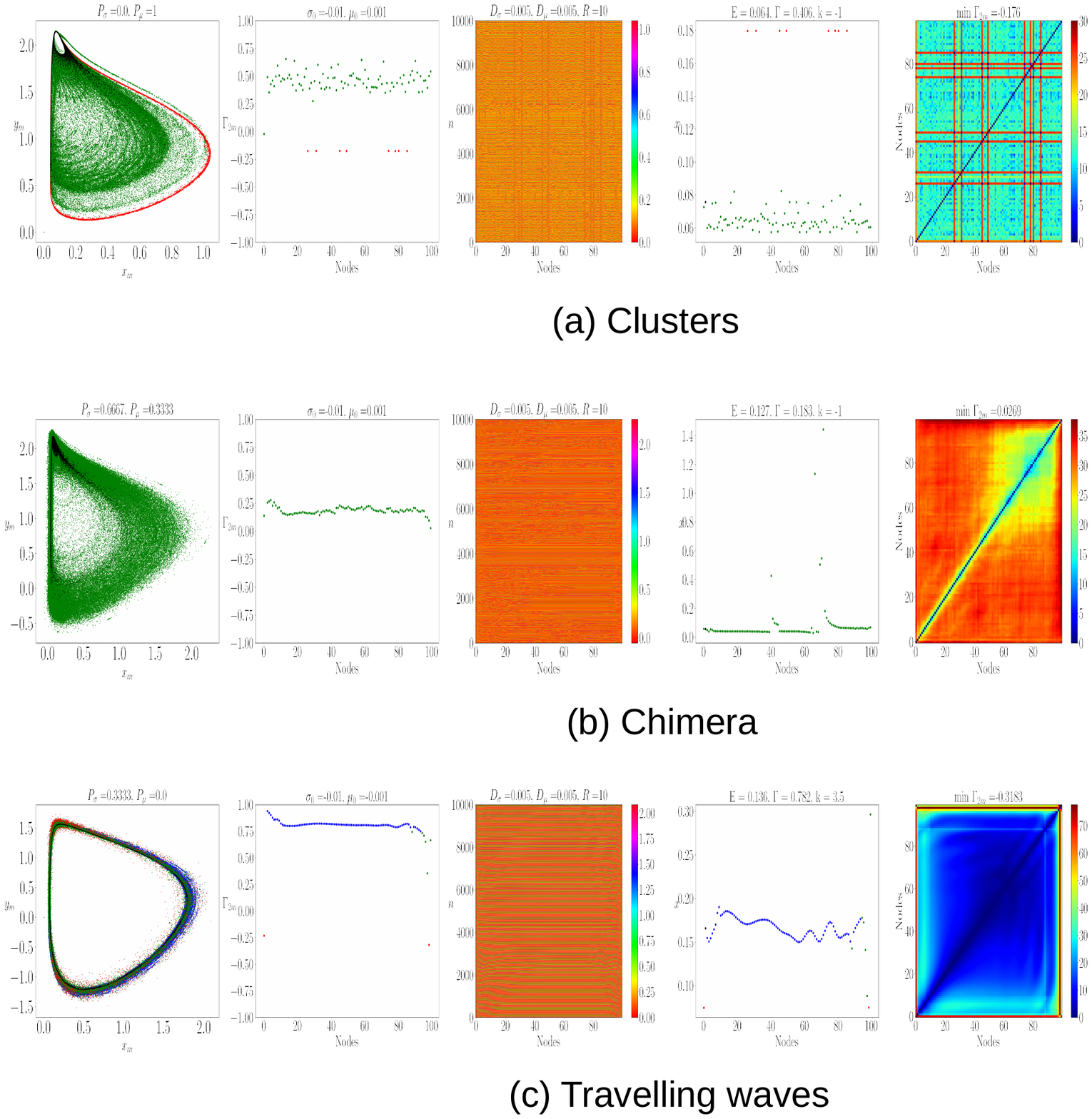}
     \caption{From left to right: phase portrait, $\Gamma_{\rm 2, m}$ plot for each node, spatiotemporal diagram, last instance plot of $x_m$ for each node, and recurrence plot. Parameters are set as (a) $\sigma_0 = -0.01, \mu_0=0.001, D_{\sigma}=0.005, D_{\mu}=0.005, P_{\sigma}=0, P_{\mu}=1,$ and $k=-1$, (b) $\sigma_0 = -0.01, \mu_0=0.001, D_{\sigma}=0.005, D_{\mu}=0.005, P_{\sigma}=0.66666, P_{\mu}=0.33333,$ and $k=-1$, and (c) $\sigma_0 =-0.01, \mu_0=-0.001, D_{\sigma}=0.005, D_{\mu}=0.005, P_{\sigma}=0.33333, P_{\mu}=0,$ and $k=3.5$. The normalized cross-correlation coefficient and the normalized syncronization error values are (a) $\Gamma = 0.406,$ $E = 0.064$, (b) $\Gamma = 0.183,$ $E = 0.127$, and (c) $\Gamma = 0.782,$ $E = 0.136$. In the first, second and fourth plots, the nodes lying in the solitary regime are denoted by red dots, the nodes with $\Gamma_{\rm 2,m} \in [-0.15, 0.75]$ are denoted by green dots and the nodes in the coherent domain are denoted by blue dots.}
 \label{fig:5plots2}
\end{figure*}

In this section, we have numerically plotted and gathered some interesting behaviors exhibited by the network under variations of different parameter values. We showcase a variety of different phase portraits, spatiotemporal patterns exhibited by the heterogeneous Chialov ring-star network in Eq.~\ref{eq:ChialvoDiscrete} and Eq.~\ref{eq:ChialvoDiscreteCentral}. For each simulation, we have shown five separate plots corresponding to different responses. In both Fig.~\ref{fig:5plots1} and Fig.~\ref{fig:5plots2}, 
\begin{enumerate}
    \item the first plot indicates the phase portrait of all the $N$ nodes with the transients removed,
    \item the second plot refers to the $\Gamma_{2,m}$ values for all the $N$ nodes,
    \item the third plot denotes the spatiotemporal dynamics of the nodes with time,
    \item the fourth plot portrays the last instance of all the units in space after a sufficient number of iterations $n$ with the transient dynamics removed, and
    \item the last plot corresponds to the recurrence plot \cite{muni2022discrete} comparing the distance between the final position of each node with all the other nodes in space after $n$ iterations.
\end{enumerate}
In the first, second, and fourth plots, the nodes lying in the solitary regime are denoted by red dots, the nodes with $\Gamma_{\rm 2,m} \in [-0.15, 0.75]$ are denoted by green dots and the nodes in the coherent domain are denoted by blue dots. Additionally, the second node is denoted by a black dot in the first (phase portrait) plot.\par

Fig.~\ref{fig:5plots1}(a) displays the solitary nodes and the coherent nodes clustered in their respective domains along with one node that rests in the region where $\Gamma_{\rm 2,m} \in [-0.15, 0.75]$. The blue nodes have cluttered with $\Gamma_{\rm 2,m} \in [0.75, 1]$ and the solitary nodes have accumulated with $\Gamma_{\rm 2,m} \sim -0.172$. We notice that the solitary nodes are distributed over the whole ensemble of the nodes. The fact that there exists almost equal number of nodes in both the clusters is evident from the tiny squares in the recurrence plot. In Fig.~\ref{fig:5plots1}(b), we see that all the nodes have been clustered and in synchrony except two, which are solitary. Synchronization is also confirmed from the very small value of the normalized synchronization-error, i.e, $E = 0.034$. Here too the solitary nodes consist of $\Gamma_{\rm 2,m} \sim -0.172$. In the corresponding recurrence plot we notice a deep blue region that covers almost the whole space, visually denoting the fact that the nodes are synchronized. Fig.~\ref{fig:5plots1}(c) has the similar dynamics as that of Fig.~\ref{fig:5plots1}(a). The solitary nodes mostly accumulate with $\Gamma_{\rm 2,m} \sim -0.167$. Note that after sufficient time iterations, the blue and the green nodes can rest together in a single cluster (See the fourth plot in Figs.~\ref{fig:5plots1}(b) and (c)). As expected, their recurrence plots also look very similar to each other. In Fig.~\ref{fig:5plots2}(a) we again see an emergence of clusters with nodes in synchrony at the solitary regime, having $\Gamma_{\rm 2,m} \sim -0.176$ and the remaining clustered around $\Gamma_{\rm 2,m} \sim 0.5$. Here in the recurrence plot, although we see squares denoting multi clusters in the dynamics, they are obviously bigger in area than the ones appearing in Fig.~\ref{fig:5plots1}(a) and Fig.~\ref{fig:5plots1}(c), due to the fact that the number of nodes in one of the clusters is much higher. The behavior in Fig.~\ref{fig:5plots2}(b) refers to the phenomenon of ``chimera" where nodes within a particular boundary (approximately in $40 \leq m \leq 45$ and in $62 \leq m \leq 75$) are completely asynchronous compared to the other nodes in the space which are convincingly synchronous, and they coexist \cite{rybalova2021interplay}. In Fig.~\ref{fig:5plots2}(c), we observe the emergence of \textit{travelling waves}. The recurrence plots in Fig.~\ref{fig:5plots2}(b) and Fig.~\ref{fig:5plots2}(c) visually support the fact that we indeed notice chimera and travelling waves, respectively.

\subsection{Time series analysis}
\begin{figure*}
\centering
   \includegraphics[scale=0.115]{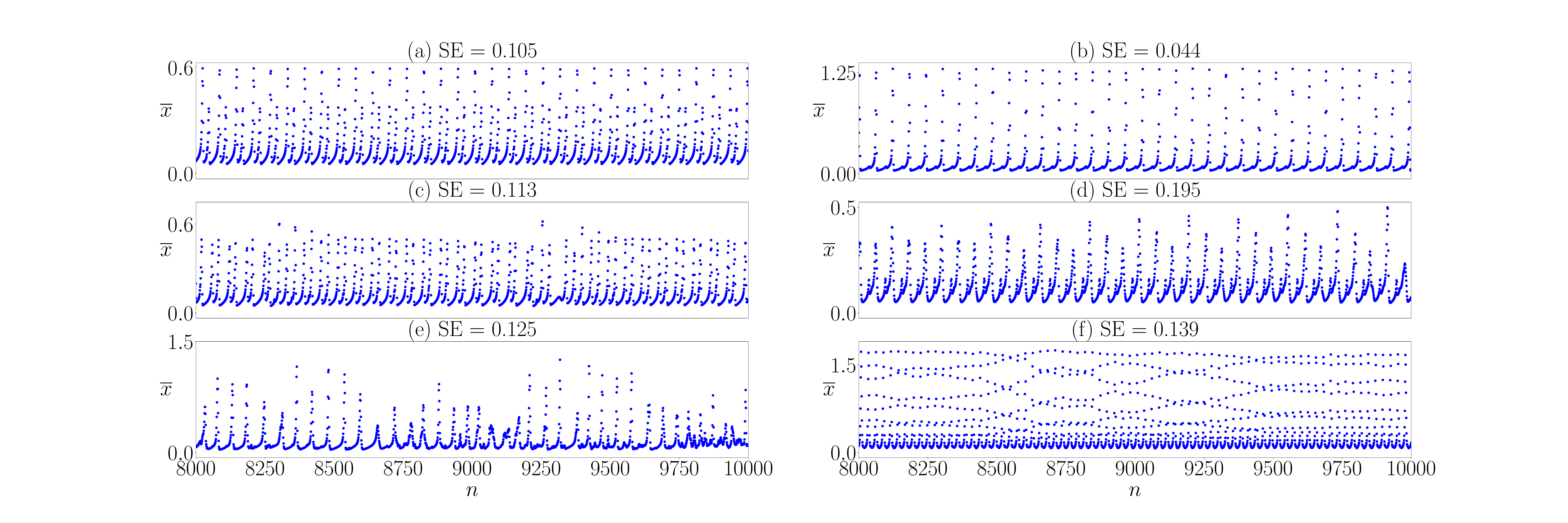}
 \caption{Time series plots along with the sample entropy values for the spatially averaged action potential obtained from Fig.~\ref{fig:5plots1} and Fig.~\ref{fig:5plots2}. Figures (a)-(c) correspond to Fig.~\ref{fig:5plots1}(a)-(c) respectively, and figures (d)-(f) correspond to Fig.~\ref{fig:5plots2}(a)-(c) respectively. SE values reported are (a) 0.105, (b) 0.044, (c) 0.113, (d) 0.195, (e) 0.125, and (f) 0.139.}
 \label{fig:TS}
\end{figure*}

Next, we perform statistical analysis on the time-series data of the action potential $x$ corresponding to the parameter combinations reported in Fig.~\ref{fig:5plots1} and 
Fig.~\ref{fig:5plots2}. As mentioned in Sec.~\ref{sec:stat_ts}, we take the spatial average of all the nodes at a particular time $n$, denoted by $\overline{x}$, and this is illustrated in Fig.~\ref{fig:TS}. The plots of the spatial average against time for the parameter combinations in Fig.~\ref{fig:5plots1}(a)-(c) and Fig.~\ref{fig:5plots2}(a)-(c) are shown in Figs.~\ref{fig:TS}(a)-(c) and Fig.~\ref{fig:TS}(d)-(f), respectively. As seen in Fig.~\ref{fig:TS}, the oscillations are irregular and do not appear to converge to any steady state of the system. To further our analysis on complex dynamics, we use \texttt{nolds}  to compute the sample entropy for each of these cases and record them on each time series plot. What we clearly observe is the presence of fairly complex behaviour in all of them, giving us an intuition about the extent of disorderedness. Visually, we also notice irregular oscillatory behaviour in the firing pattern of $\overline{x}$, not converging to any stable steady state.  Note that Fig.~\ref{fig:TS}(b) which corresponds to the mostly synchronized case ($E = 0.034$) Fig.~\ref{fig:5plots1}(b) has the lowest value of sample entropy (SE=$0.044$) and thus the lowest disorderedness as compared to the other five cases. Furthermore, out of the above six cases, the highest value of synchronization error is observed in Fig.~\ref{fig:5plots2}(c), having $E=0.136$. The corresponding time series Fig.~\ref{fig:TS}(a) too has a very high value of sample entropy, SE=$0.139$ (second highest). Statistically speaking, it can be expected that an increase in asynchrony leads to an increase in the complexity of the system dynamics. From the color plots and their corresponding $E$ vs. SE plots in the next section, we can infer the above phenomenon.

\subsection{Two dimensional color coded plots}

\subsubsection{Parameter space defined by $(\sigma_0, \mu_0)$}
\begin{figure}
\centering
    \centering    \includegraphics[scale=0.07]{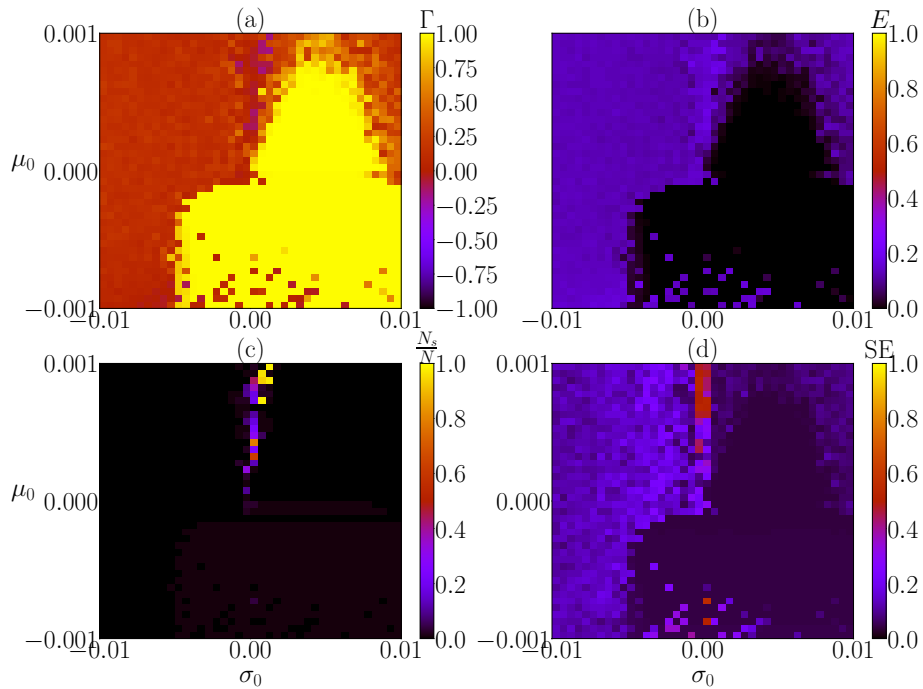}
 \caption{Collection of two-dimensional color plots of (a) $\Gamma,$ (b) $E,$ (c) $\frac{N_s}{N},$ and (d) SE in the parameter space defined by the mean coupling strengths $(\sigma_0, \mu_0)$. The parameter space is a $40 \times 40$ grid with $k=-1, P_{\mu}=1, P_{\sigma} \sim 0.666, D_{\mu}=D_{\sigma}=0.005$. An almost definitive bifurcation boundary within which most of the nodes lie in the coherent regime, is synchronous and has low sample entropy, is observed. Solitary nodes appear when the mean coupling strengths are $\sigma_0 \sim 0$ and $\mu_0>0$, corresponding the region where the nodes are asynchronous and have high sample entropy.}
 \label{fig:sigma_0vsmu_0colorplot}
\end{figure}

\begin{figure}
\centering
  \includegraphics[scale=0.05]{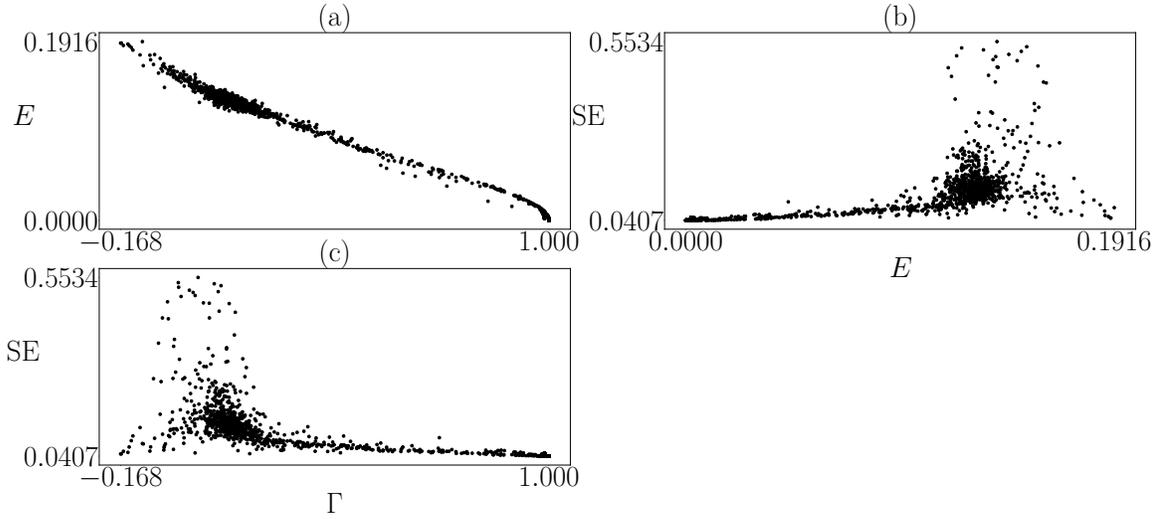}
 \caption{Comparison plots for the various measures: (a) $E$ vs. $\Gamma$, (b) SE vs. $E$, and (c) SE vs. $\Gamma$, collected from Fig.~\ref{fig:sigma_0vsmu_0colorplot}. Figures (a) and (c) show an inverse trend whereas figure (b) shows a proportional trend.}
 \label{fig:sigma_0vsmu_0comp}
\end{figure}

Fig.~\ref{fig:sigma_0vsmu_0colorplot} is a collection of two-dimensional color plots of $\Gamma, E, \frac{N_s}{N}$, and SE in the parameter space defined by $(\sigma_0, \mu_0)$. The parameter space is a $40 \times 40$ grid with $k=-1, P_{\mu}=1, P_{\sigma} \sim 0.666, D_{\mu}=D_{\sigma}=0.005$. From the plot that depicts the normalized cross correlation coefficient (Fig.~\ref{fig:sigma_0vsmu_0colorplot}(a)), we observe an almost definitive bifurcation boundary within which most of the nodes lie in the coherent regime, i.e, $\Gamma \sim 1$ and outside which the nodes behave incoherently, i.e, $\Gamma < 0.75$. A similar kind of bifurcation boundary is observed in the normalized syncronization error color plot (Fig.~\ref{fig:sigma_0vsmu_0colorplot}(b)). The region where $\Gamma \sim 1$, the nodes tend to adopt a complete synchronization behavior, making $E \sim 0$. Moving on to the sample entropy plot (Fig.~\ref{fig:sigma_0vsmu_0colorplot}(d)), interestingly we once again notice an almost similar bifurcation boundary between no complexity (SE$\sim 0$) and onset of complexity (SE $>0$). Whenever, $\Gamma \sim 1$ and $E \sim 0$ then SE$\sim 0$ too. Now in the parameter region where $\Gamma<1$, we notice a transition in the behavior of the nodes from complete synchrony to a degree of asynchrony, such that the normalized synchronization rate lies in the range $0<E<0.3$ approximately. In an analogous mechanism the sample entropy increases too when $\Gamma < 1$ pertaining to a more disordered system dynamics. When $\Gamma < 1$, and the mean coupling strength $\sigma_0 \sim 0$, there appears a tinge of violet region in the color plot, implying $\Gamma \in (-0.15, -0.38)$, for which there appear solitary nodes (as evident from the $\frac{N_s}{N}$ color plot, Fig.~\ref{fig:sigma_0vsmu_0colorplot}(c)) and a further peak in the value of sample entropy ($0.3< {\rm SE} < 0.75$ approximately). Otherwise throughout the parameter space, there exit no solitary nodes. In Fig.~\ref{fig:sigma_0vsmu_0comp} we have collected the values of $E, \Gamma$, and SE and plotted them against each other. Notice in Fig.~\ref{fig:sigma_0vsmu_0comp}(a) that it gives a clear inversely proportional trend for $E$ and $\Gamma$. In Fig.~\ref{fig:sigma_0vsmu_0comp}(b), we notice that with increase in $E$, the sample entropy shows a fairly increasing trend, whereas in Fig.~\ref{fig:sigma_0vsmu_0comp}(c), with increase in $\Gamma$, the sample entropy decreases as expected.

\subsubsection{Parameter space defined by $(\sigma_0, D_{\sigma})$}
\begin{figure}[!htbp]
\centering   \includegraphics[scale=0.07]{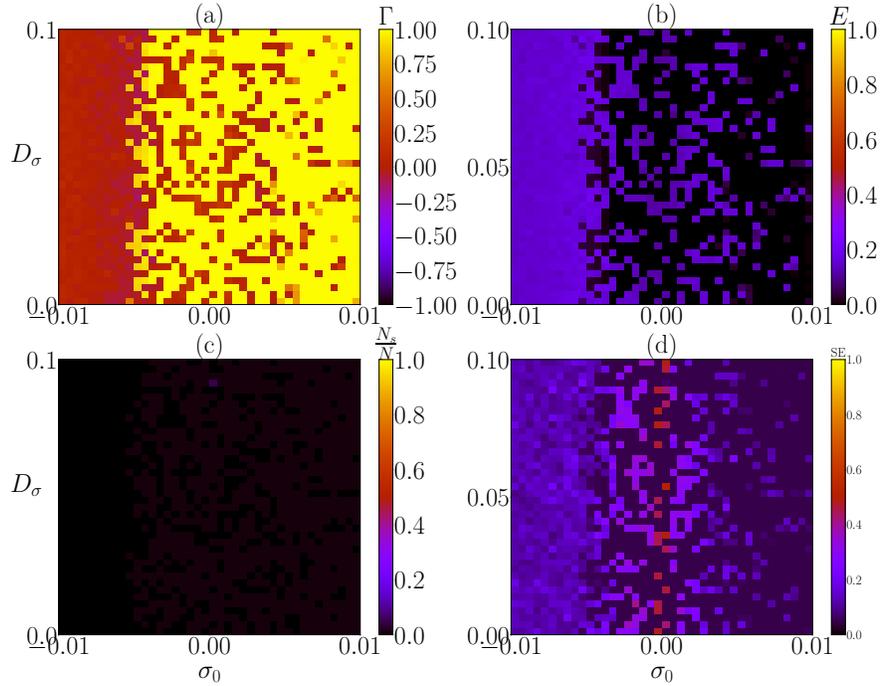}
 \caption{Collection of two-dimensional color plots of (a) $\Gamma,$ (b) $E,$ (c) $\frac{N_s}{N},$ and (d) SE in the parameter space defined by the mean coupling strength and the noise modulation $(\sigma_0, D_{\sigma})$. The parameter space is a $40 \times 40$ grid with $\mu_0=-0.001, P_{\mu}=1, P_{\sigma} \sim 0.666, D_{\mu}=0.005$ and $k=-1$. For $-0.01<\sigma_0<-0.005$, nodes have $\Gamma \sim (0, 0.5)$ and are asynchronous. When $\sigma_0>-0.005$, the nodes mostly tend towards synchronization. Sample entropy is the highest at $\sigma_0 \sim 0$. Solitary nodes rarely appear.}
 \label{fig:sigma_0vsD_sigmacolorplot}
\end{figure}

\begin{figure}
\centering
  \includegraphics[scale=0.05]{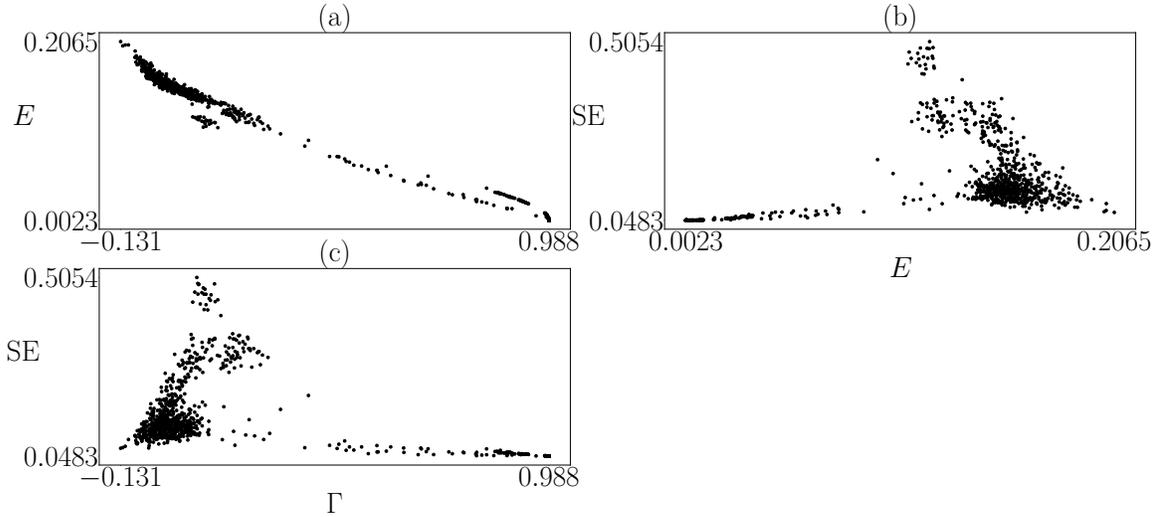}
 \caption{Comparison plots for the various measures: (a) $E$ vs. $\Gamma$, (b) SE vs. $E$, and (c) SE vs. $\Gamma$, collected from Fig.~\ref{fig:sigma_0vsD_sigmacolorplot}. Figures (a) and (c) show an inverse trend whereas figure (b) shows a proportional trend.}
 \label{fig:sigma_0vsD_sigmacomp}
\end{figure}

Fig.~\ref{fig:sigma_0vsD_sigmacolorplot} is a collection of two-dimensional color plots in the parameter space defined by $(\sigma_0, D_{\sigma})$ with $\mu_0=-0.001, P_{\mu}=1, P_{\sigma} \sim 0.666, D_{\mu}=0.005$ and $k=-1$. Note that for all values of $D_{\sigma}$, when $-0.01<\sigma_0<-0.005$, we see that the nodes prefer to remain distributed throughout the space, with $\Gamma \sim (0, 0.5)$ (Fig.~\ref{fig:sigma_0vsD_sigmacolorplot}(a)) and behave in an asynchronous manner as evident from the violet region in the $E$ color plot (Fig.~\ref{fig:sigma_0vsD_sigmacolorplot}(b)) denoting higher asynchronity. As soon as $\sigma_0>-0.005$, the nodes mostly tend towards synchronization in the coherent region as can be seen from the deep yellow and black regimes respectively in the $\Gamma$ and $E$ color plots respectively, although there might be instances when the nodes do not settle to synchronity in the coherent domain. Interestingly, there exists a very few to no solitary nodes as can be seen from the $\frac{N_s}{N}$ color plot (Fig.~\ref{fig:sigma_0vsD_sigmacolorplot}(c)). Now looking into the SE plot (Fig.~\ref{fig:sigma_0vsD_sigmacolorplot}(d)), note that again there is a fair proportional relationship between sample entropy and asynchronicity, with the maximum complex behaviour appearing at $\sigma_0 \sim 0$. It can be seen that with increase in the mean coupling strength $\sigma_0$, from $-0.01$ to $0.01$, sample entropy starts increasing, reaches a peak at around $\sigma_0 \sim 0$ and then the value drops. In Fig.~\ref{fig:sigma_0vsD_sigmacomp} we have collected the values of $E, \Gamma$, and SE and plotted them against each other. Notice in Fig.~\ref{fig:sigma_0vsD_sigmacomp}(a) that it clearly indicates an inversely proportional trend for $E$ and $\Gamma$. In Fig.~\ref{fig:sigma_0vsD_sigmacomp}(b), we notice that with increase in $E$, the sample entropy shows a fairly increasing trend, whereas in Fig.~\ref{fig:sigma_0vsD_sigmacomp}(c), with increase in $\Gamma$, the sample entropy decreases as expected. 

\subsubsection{Parameter space defined by $(\mu_0, D_{\mu})$}
\begin{figure}
\centering
   \includegraphics[scale=0.07]{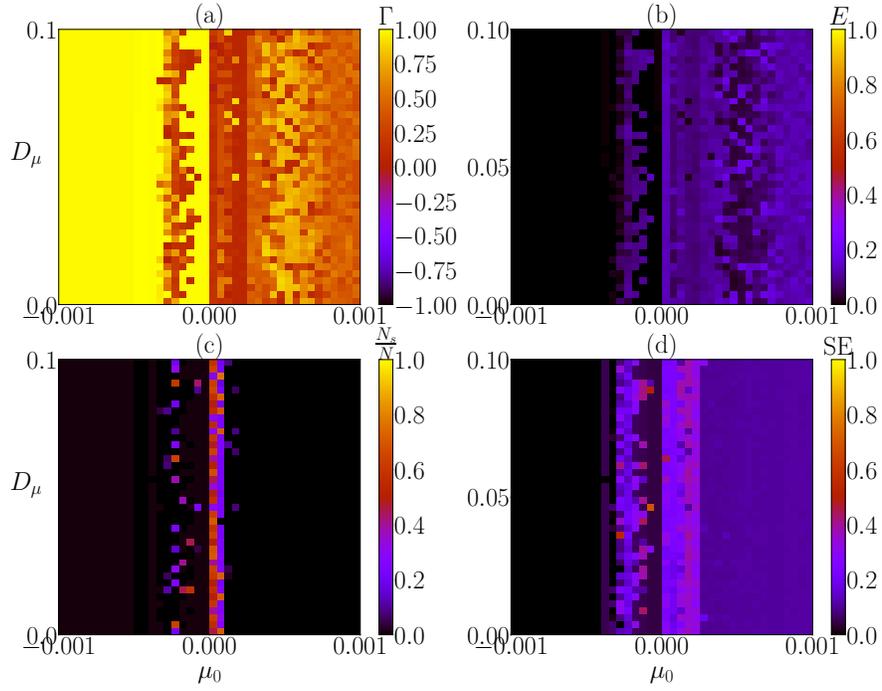}
 \caption{Collection of two-dimensional color plots of (a) $\Gamma,$ (b) $E,$ (c) $\frac{N_s}{N},$ and (d) SE in the parameter space defined by the mean coupling strength and the noise modulation $(\mu_0, D_{\mu})$. The parameter space is a $40 \times 40$ grid with $\sigma_0=0, P_{\mu}=1, P_{\sigma} \sim 0.666, D_{\sigma}=0.005$ and $k=-1$. For all values of $D_{\mu} > 0$, probability of all the nodes having $\Gamma<0.75$ is very high. Solitary nodes appear where  $-0.38<\Gamma<-0.15$, i.e, at around $\mu_0 \sim 0$. Asynchronity and high disorderedness are observed for $\Gamma>0.75$.}
 \label{fig:mu_0vsD_mucolorplot}
\end{figure}

\begin{figure}
\centering
  \includegraphics[scale=0.05]{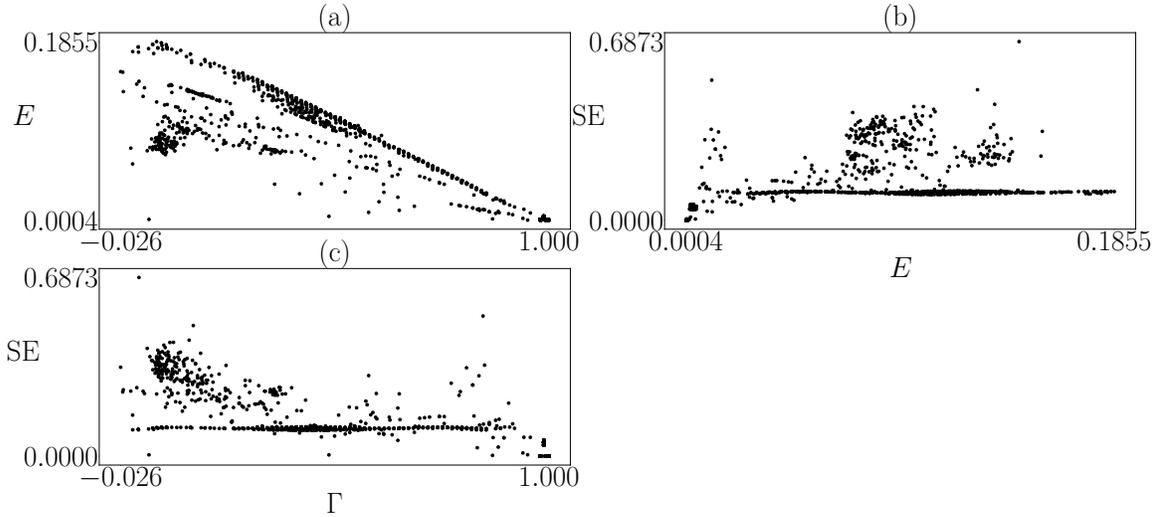}

 \caption{Comparison plots for the various measures: (a) $E$ vs. $\Gamma$, (b) SE vs. $E$, and (c) SE vs. $\Gamma$, collected from Fig.~\ref{fig:mu_0vsD_mucolorplot}. Figures (a) and (c) show an inverse trend whereas figure (b) shows a proportional trend.}
 \label{fig:mu_0vsD_mucomp}
\end{figure}

Fig.~\ref{fig:mu_0vsD_mucolorplot} is a collection of two-dimensional color plots in the parameter space defined by $(\mu_0, D_{\mu})$ with $\sigma_0=0, P_{\mu}=1, P_{\sigma} \sim 0.666, D_{\sigma}=0.005$ and $k=-1$. We note that for all values of $D_{\mu} > 0$, probability of all the nodes having $\Gamma<0.75$ is very high, as depicted by reddish-yellow region in Fig.~\ref{fig:mu_0vsD_mucolorplot} (a), with $-0.38<\Gamma<-0.15$ at around $\mu_0 \sim 0$, as can be seen from the emergence of solitary nodes in Fig.~\ref{fig:mu_0vsD_mucolorplot} (c) in this region. Otherwise there occurs no solitary nodes for $\mu_0>0$. For $\mu_0 < 0$, mostly there exists nodes with $\Gamma \sim 1$ denoting nodes in the coherent domain and synchronized (Fig.~\ref{fig:mu_0vsD_mucolorplot} (b)), except a certain region of $\mu_0<0$, where again solitary nodes emerge. Nodes are synchronized where $\Gamma>0.75$. From the sample entropy plot (Fig.~\ref{fig:mu_0vsD_mucolorplot} (d)) we observe that for $\mu_0>0$, there occurs high disorderdness, i.e, for the paramter combinations where $\gamma < 0.75$ and nodes are asynchronous. Again, like $\Gamma$ and $E$ plot, there exists a region $\mu <0$, where the sample entropy is $>0$, overlapping with the region where solitary nodes emerge. In Fig.~\ref{fig:mu_0vsD_mucomp} we have collected the values of $E, \Gamma$, and SE and plotted them against each other. Notice in Fig.~\ref{fig:mu_0vsD_mucomp}(a) that it clearly indicates an inversely proportional trend for $E$ and $\Gamma$. In Fig.~\ref{fig:mu_0vsD_mucomp}(b), we notice that with increase in $E$, the sample entropy shows a fairly increasing trend, whereas in Fig.~\ref{fig:mu_0vsD_mucomp}(c), with increase in $\Gamma$, the sample entropy decreases as expected. 

\subsubsection{Parameter space defined by $(P_{\mu}, P_{\sigma})$}
\begin{figure}[!htbp]
\centering
   \includegraphics[scale=0.07]{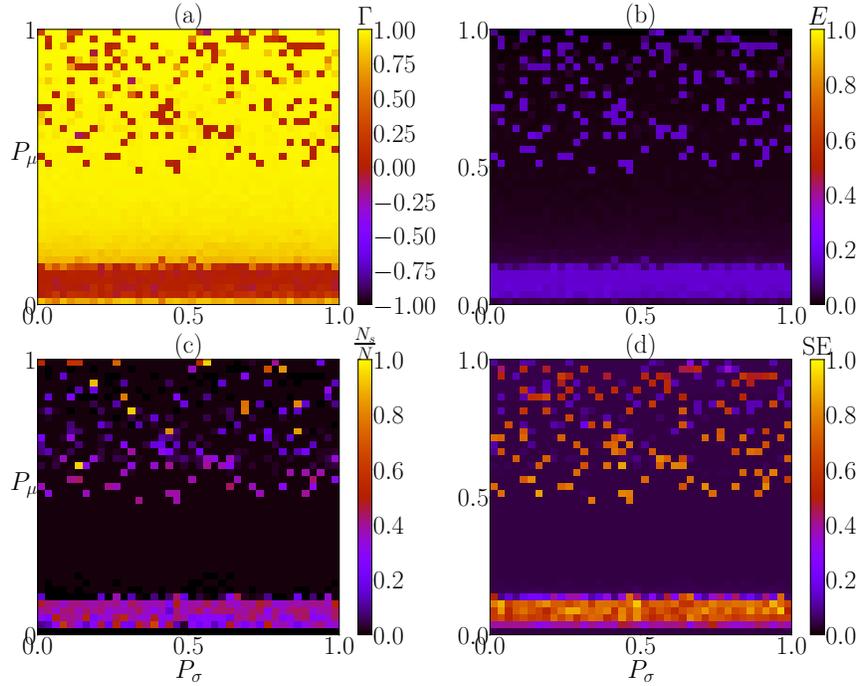}
 \caption{Collection of two-dimensional color plots of (a) $\Gamma,$ (b) $E,$ (c) $\frac{N_s}{N},$ and (d) SE in the parameter space defined by the coupling probabilities $(P_{\sigma}, P_{\mu})$. The parameter space is a $40 \times 40$ grid with $\sigma_0=0, \mu_0 = -0.001, D_{\sigma}=D_{\mu}=0.005$ and $k=-1$. For values close to $P_{\mu}=0$, for all $P_{\sigma}$, there exist regions where $\Gamma<0.75$, and the nodes behave asynchronously, having very high sample entropy value and are solitary. As soon as $P_{\mu}>0.5$, the region is randomly distributed between both types of extreme behaviors.}
 \label{fig:p_sigmavsp_mucolorplot}
\end{figure}

\begin{figure}
\centering
 \includegraphics[scale=0.05]{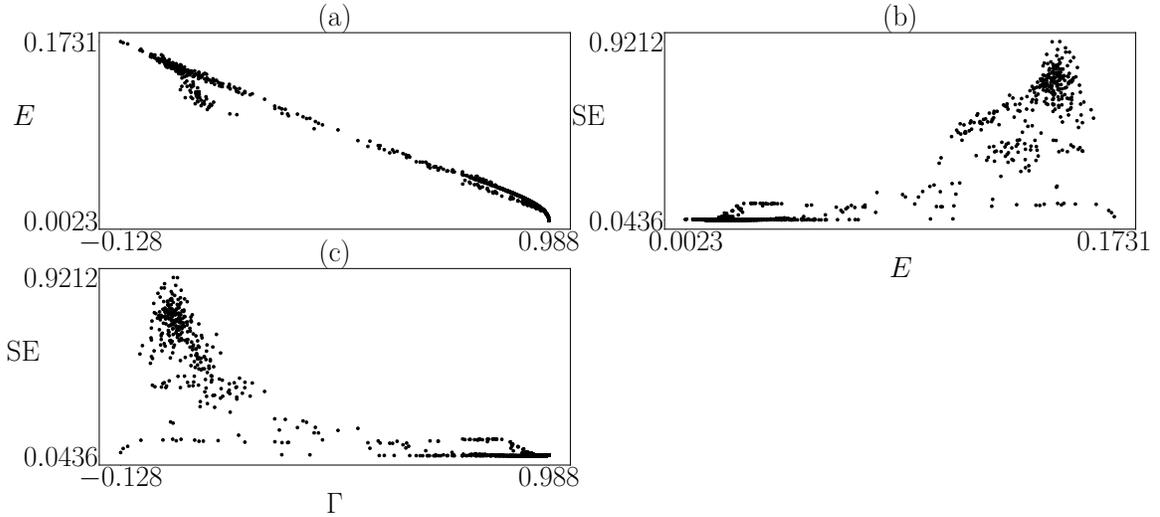}
 \caption{Comparison plots for the various measures: (a) $E$ vs. $\Gamma$, (b) SE vs. $E$, and (c) SE vs. $\Gamma$, collected from Fig.~\ref{fig:p_sigmavsp_mucolorplot}. Figures (a) and (c) show an inverse trend whereas figure (b) shows a proportional trend.}
 \label{fig:p_sigmavsp_mucomp}
\end{figure}

Fig.~\ref{fig:p_sigmavsp_mucolorplot} is a collection of two-dimensional color plots in the parameter space defined by $(P_{\mu}, P_{\sigma})$ with $\sigma_0=0, \mu_0 = -0.001, D_{\sigma}=D_{\mu}=0.005$ and $k=-1$. From Fig.~\ref{fig:p_sigmavsp_mucolorplot}(a), we see that for values close to $P_{\mu}=0$, for all $P_{\sigma}$, there exist regions where $\Gamma<0.75$, where nodes behave asynchronously (Fig.~\ref{fig:p_sigmavsp_mucolorplot}(b)), having very high sample entropy value (Fig.~\ref{fig:p_sigmavsp_mucolorplot}(d)) and are solitary (Fig.~\ref{fig:p_sigmavsp_mucolorplot}(c)). Above this region of $P_{\mu}$ where $P_{\mu}<0.5$, mostly the nodes are synchronized in the coherent domain having very small value of sample entropy. As soon as $P_{\mu}>0.5$, the region is randomly distributed between both types of extreme behaviors, denoted by the contrasting colored boxes in all the four color coded plots. In Fig.~\ref{fig:p_sigmavsp_mucomp} we have collected the values of $E, \Gamma$, and SE and plotted them against each other. Notice in Fig.~\ref{fig:p_sigmavsp_mucomp}(a) that it clearly indicates an inversely proportional trend for $E$ and $\Gamma$. In Fig.~\ref{fig:p_sigmavsp_mucomp}(b), we notice that with increase in $E$, the sample entropy shows a fairly increasing trend, whereas in Fig.~\ref{fig:p_sigmavsp_mucomp}(c), with increase in $\Gamma$, the sample entropy decreases as expected. 

\subsubsection{Parameter space defined by $(\sigma_0, k)$}
\begin{figure}[!htbp]
\centering
 \includegraphics[scale=0.07]{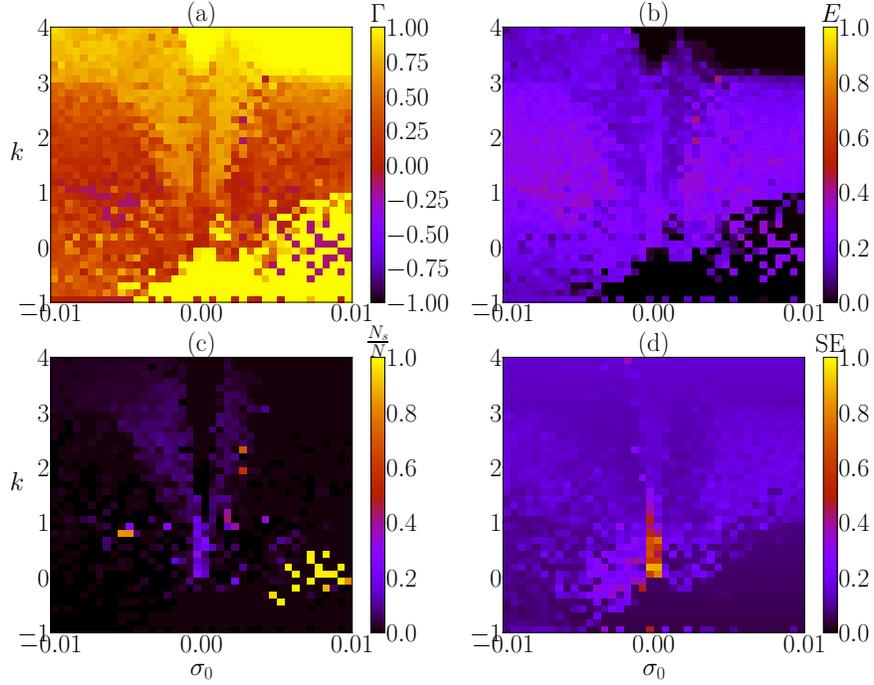}
 \caption{Collection of two-dimensional color plots of (a) $\Gamma,$ (b) $E,$ (c) $\frac{N_s}{N},$ and (d) SE in the parameter space defined by the mean coupling strength and the electromagnetic flux coupling $(\sigma_0, k)$. The parameter space is a $40 \times 40$ grid with $\mu_0=-0.001, P_{\mu}=1, P_{\sigma} \sim 0.666, D_{\mu}=D_{\sigma}=0.005$. The nodes mostly behave in a fairly asynchronous manner having $\Gamma $values that lie in $(-0.38, 0.75)$. There exist two moderately distinct regions in the where the nodes all lie in the coherent region and they are almost completely synchronous. The sample entropy is the highest at $\sigma_0 \sim 0$ and $k \in (0, 1.5)$. there exist some regions where all the nodes are solitary within $k \in (-0.5, 0.5), \sigma_0 \in (0.005, 0.01)$}
 \label{fig:sigma_0vskcolorplot}
\end{figure}

\begin{figure}
\centering
 \includegraphics[scale=0.05]{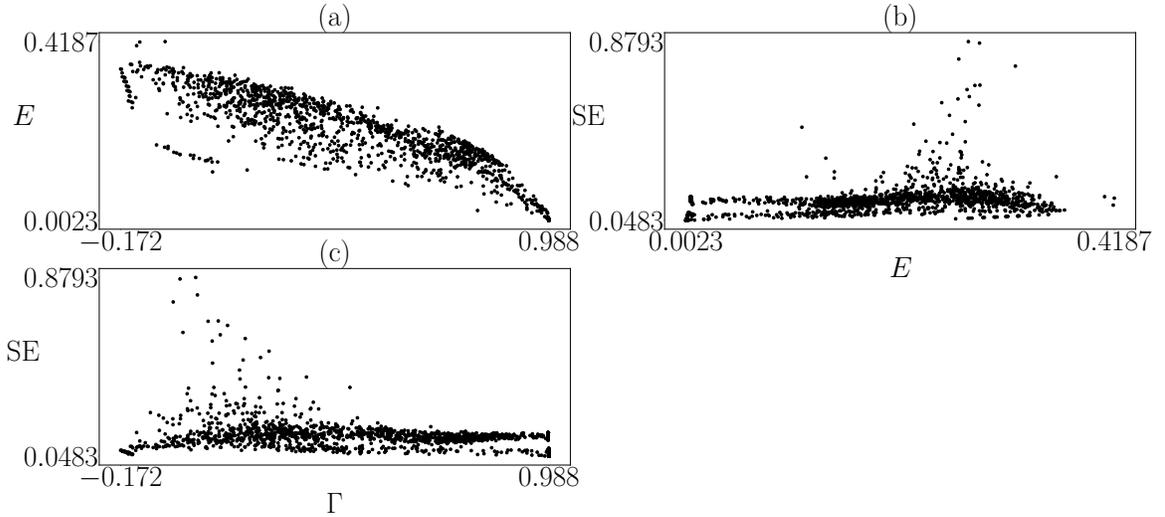}
 \caption{Comparison plots for the various measures: (a) $E$ vs. $\Gamma$, (b) SE vs. $E$, and (c) SE vs. $\Gamma$, collected from Fig.~\ref{fig:sigma_0vskcolorplot}. Figures (a) and (c) show an inverse trend whereas figure (b) shows a proportional trend.}
 \label{fig:sigma_0vskcomp}
\end{figure}

Fig~\ref{fig:sigma_0vskcolorplot} is a collection of two-dimensional color plots in the parameter space defined by $(\sigma_0, k)$ with $\mu_0=-0.001, P_{\mu}=1, P_{\sigma} \sim 0.666, D_{\mu}=D_{\sigma}=0.005$. We observe that the $\Gamma$ colored plot (Fig~\ref{fig:sigma_0vskcolorplot}(a)) is mostly dominated by values that lie in $(-0.38, 0.75)$, for which the nodes behave in a fairly asynchronous manner as depicted from the corresponding synchronization error plot $E$ (Fig~\ref{fig:sigma_0vskcolorplot}(b)). There exist two moderately distinct regions in the $\Gamma$ plot where the nodes all lie in the coherent region and they are almost completely synchronous (see the deep yellow regions and the deep black regions in the $\Gamma$ and $E$ plots respectively). Comparing the synchronization error plot and the sample entropy plot, we observe again a fairly increasing relationship between the two measures. The sample entropy (Fig~\ref{fig:sigma_0vskcolorplot}(d)) has the maximum value at the region $\sigma_0 \sim 0$ and $k \in (0, 1.5)$. Looking into the plot depicting the normalized number of solitary nodes (Fig~\ref{fig:sigma_0vskcolorplot}(c)), we detect a space where almost all the nodes are solitary in the region $k \in (-0.5, 0.5), \sigma_0 \in (0.005, 0.01)$. In Fig.~\ref{fig:sigma_0vskcomp} we have collected the values of $E, \Gamma$, and SE and plotted them against each other. Notice in Fig.~\ref{fig:sigma_0vskcomp}(a) that it clearly indicates an inversely proportional trend for $E$ and $\Gamma$. In Fig.~\ref{fig:sigma_0vskcomp}(b), we notice that with increase an in $E$, the sample entropy shows a fairly increasing trend, whereas in Fig.~\ref{fig:sigma_0vskcomp}(c), with an increase in $\Gamma$, the sample entropy decreases as expected.

\subsubsection{Parameter space defined by $(\mu_0, k)$}
\begin{figure}[!htbp]
\centering
 \includegraphics[scale=0.07]{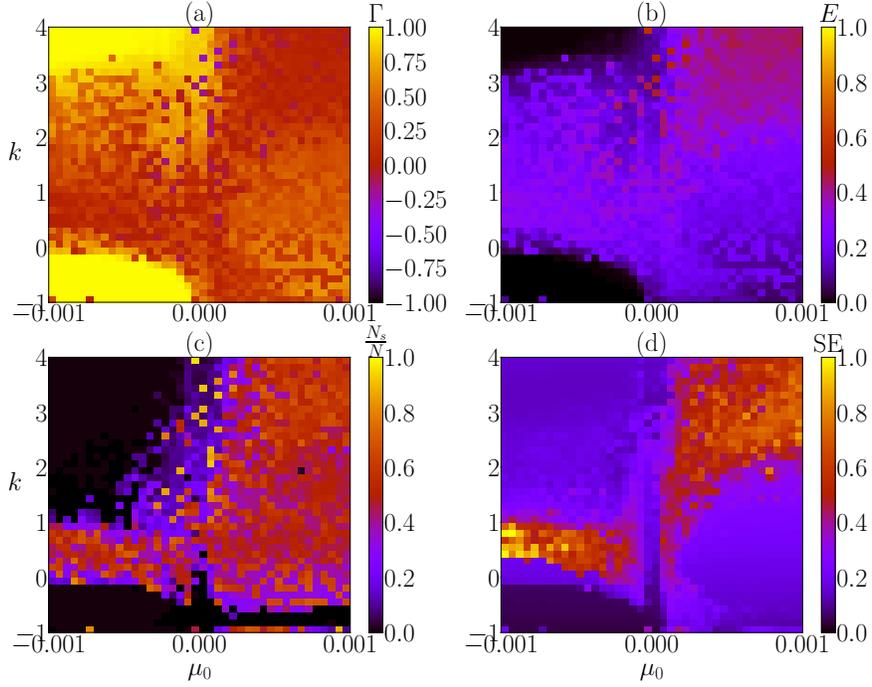}
 \caption{Collection of two-dimensional color plots of (a) $\Gamma,$ (b) $E,$ (c) $\frac{N_s}{N},$ and (d) SE in the parameter space defined by the mean coupling strength and the electromagnetic flux coupling $(\sigma_0, k)$. The parameter space is a $40 \times 40$ grid with $\sigma_0=0, P_{\mu}=1, P_{\sigma} \sim 0.666, D_{\mu}=D_{\sigma}=0.005$. For $\mu_0>0$, $\Gamma$ lies in the range  $(-0.38, 0.75)$ and the nodes behave in asynchronous fashion. $E$ and Sample entropy show high values in the region $\mu_0>0, k>1.5$. For $mu>0$, there is a high possibility of occurrence of solitary nodes, spanning over $k$ with $0.2<\frac{N_s}{N} <0.8$. When $\mu_0<0$, solitary nodes appear with $0.2<\frac{N_s}{N}<0.8$ within $0.5<k<1.5$, where there is a peak in the value of sample entropy.}
 \label{fig:mu_0vskcolorplot}
\end{figure}

\begin{figure}
\centering 
\includegraphics[scale=0.05]{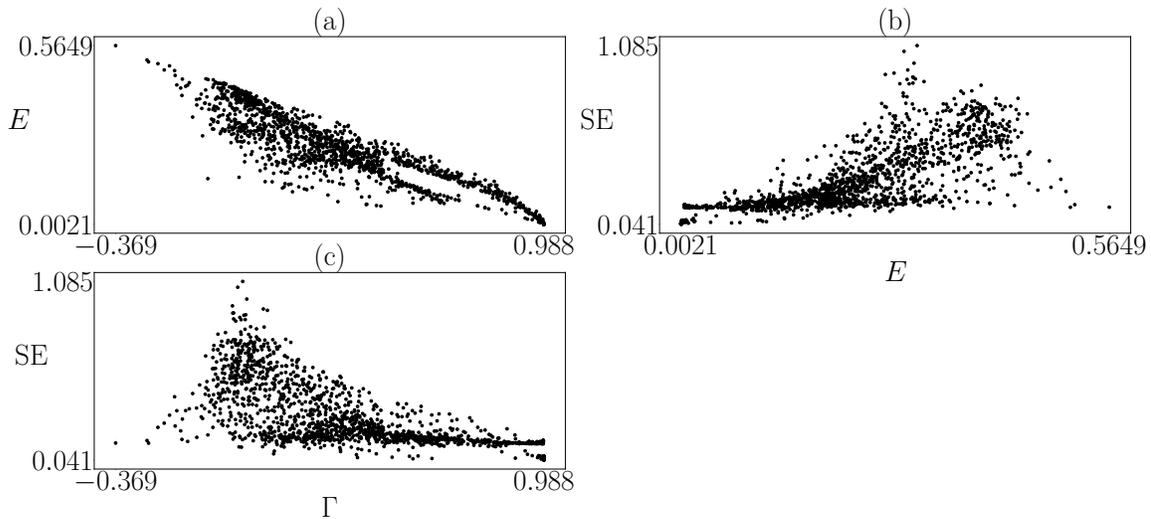}
 \caption{Comparison plots for the various measures: (a) $E$ vs. $\Gamma$, (b) SE vs. $E$, and (c) SE vs. $\Gamma$, collected from Fig.~\ref{fig:mu_0vskcolorplot}. Figures (a) and (c) show an inverse trend whereas figure (b) shows a proportional trend.}
 \label{fig:mu_0vskcomp}
\end{figure}

Fig.~\ref{fig:mu_0vskcolorplot}  is a collection of two-dimensional color plots in the parameter space defined by $(\mu_0, k)$ with $\sigma_0=0, P_{\mu}=1, P_{\sigma} \sim 0.666, D_{\mu}=D_{\sigma}=0.005$. It can be seen that for $\mu_0>0$, $\Gamma$ lies in the range  $(-0.38, 0.75)$ (Fig.~\ref{fig:mu_0vskcolorplot}(a)) and the nodes behave in asynchronous fashion as evident from a larger value in the averaged synchronization error (Fig.~\ref{fig:mu_0vskcolorplot}(b)). The value of $E$ is maximum in the region $\mu_0>0, k>1.5$, where the sample entropy (Fig.~\ref{fig:mu_0vskcolorplot}(d)) reaches very high values too. When $\mu < 0$, mostly all the nodes seem to be synchronous ($E \sim 0$) and lie in the coherent region ($\Gamma \sim 1$) for $-1<k<0$ and $4>k>2$ as can be interpreted from the $\Gamma$ and $E$ color plots. Looking into the SE plot again, for $\mu_0<0$ and $-1<k<0$, there exists a region where the sample entropy is very close to $0$. For $\mu>0$, there is a high probability of occurrence of solitary nodes (Fig.~\ref{fig:mu_0vskcolorplot}(c)) spanning over $k$ with $0.2<\frac{N_s}{N} <0.8$. When $\mu_0<0$, solitary nodes appear with $0.2<\frac{N_s}{N}<0.8$ within $0.5<k<1.5$. We also notice an upsurge in the value of the sample entropy at this parameter region. In Fig.~\ref{fig:mu_0vskcomp} we have collected the values of $E, \Gamma$, and SE and plotted them against each other. Notice in Fig.~\ref{fig:mu_0vskcomp}(a) that it clearly indicates an inversely proportional trend for $E$ and $\Gamma$. In Fig.~\ref{fig:mu_0vskcomp}(b), we notice that with increase in $E$, the sample entropy shows a fairly increasing trend, whereas in Fig.~\ref{fig:mu_0vskcomp}(c), with increase in $\Gamma$, the sample entropy decreases as expected.

\subsection{Bifurcation diagrams for synchronization}

\begin{figure*}
\begin{tabular}{cc}
  \includegraphics[scale=0.15]{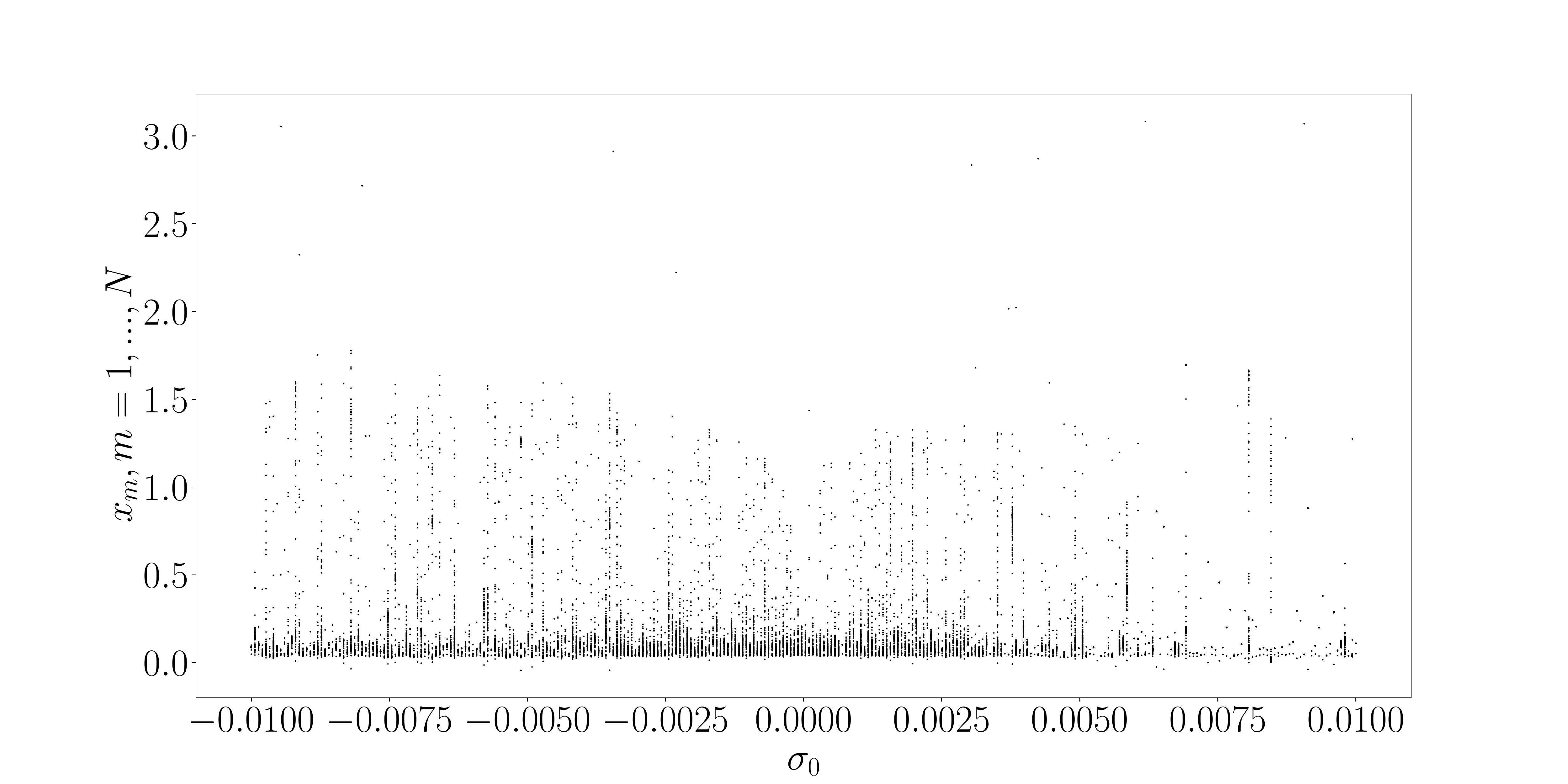} &   \includegraphics[scale=0.15]{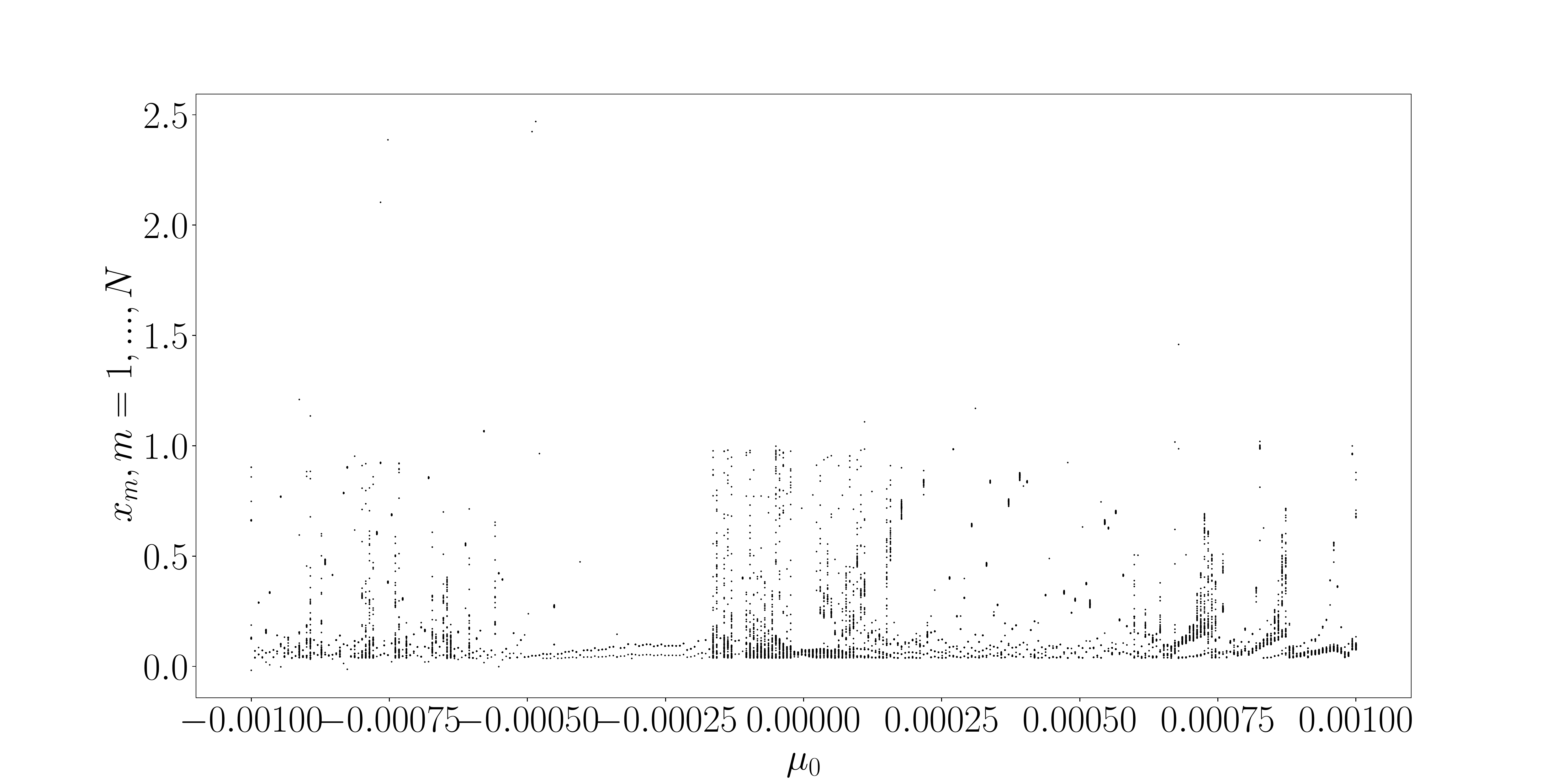} \\
(a) Varying $\sigma_0$ & (b)Varying $\mu_0$ \\[6pt]
 \includegraphics[scale=0.15]{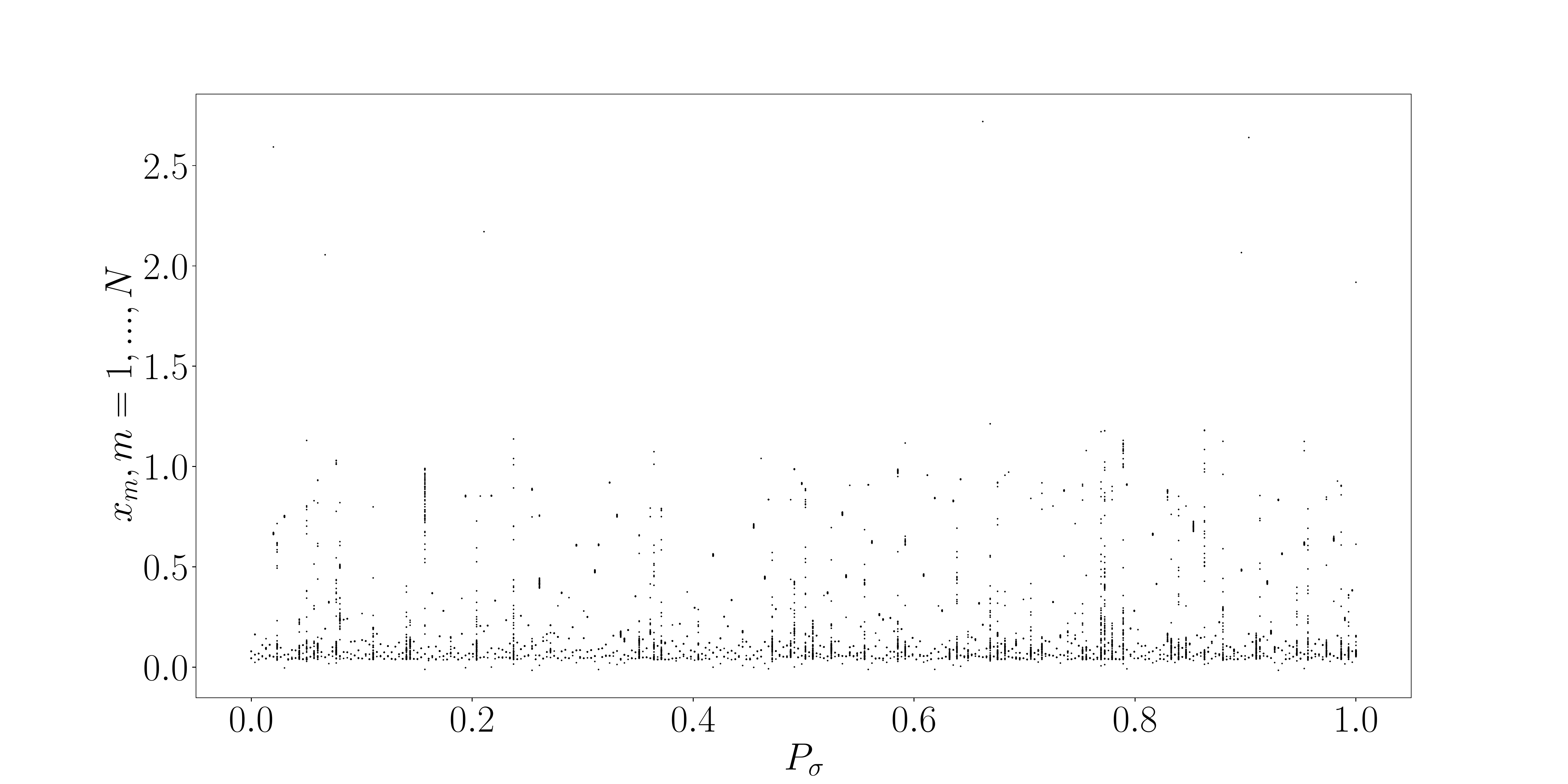} &   \includegraphics[scale=0.15]{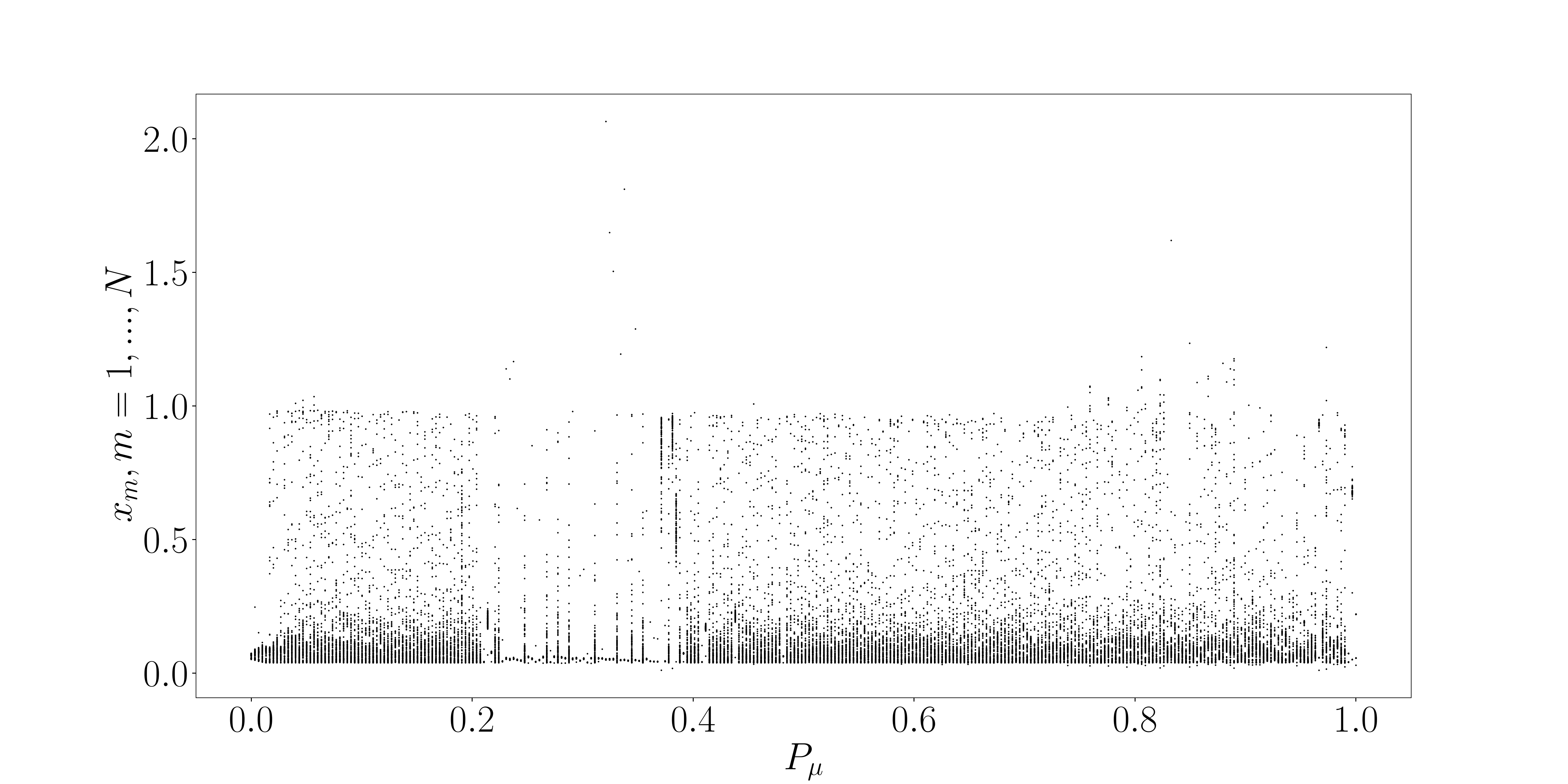} \\
(c) Varying $P_{\sigma}$ & (d) Varying $P_{\mu}$ \\[6pt]
\multicolumn{2}{c}{\includegraphics[scale=0.15]{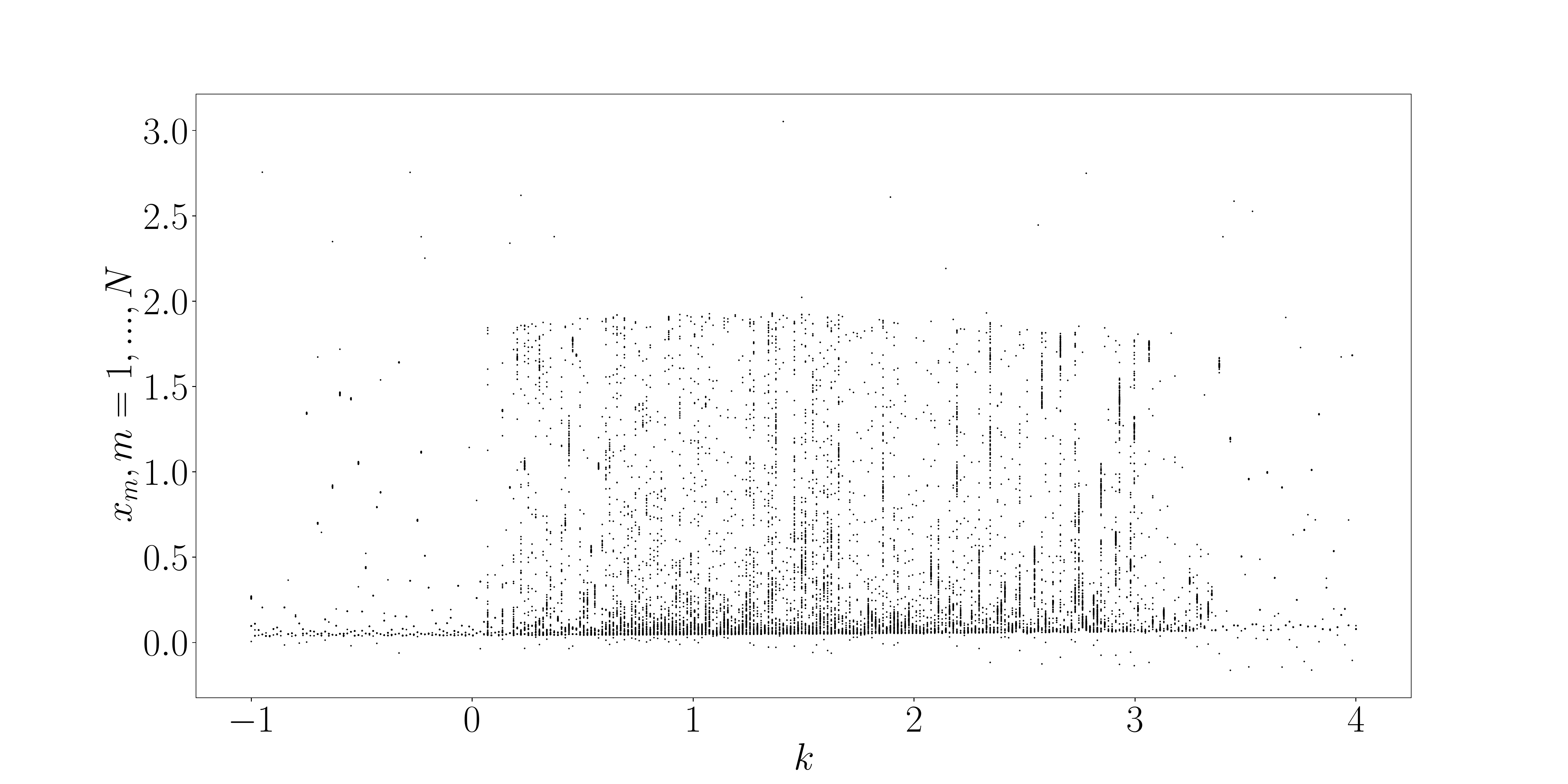} }\\
\multicolumn{2}{c}{(e) Varying $k$}
\end{tabular}
\caption{Bifurcation diagrams for the last instances of $x_m$ against the parameters  (a) $\sigma_0$, (b) $\mu_0$, (c) $P_{\sigma}$, (d) $P_{\mu}$, and (e) $k$ showing windows of synchronization and asynchronization. Other parameters are set as (a) $\mu_0=-0.001, D_{\sigma}=0.005, D_{\mu}=0.005, P_{\sigma}=0.66666, P_{\mu}=1, k=-1$, (b) $\sigma_0=0, D_{\sigma}=0.005, D_{\mu}=0.005, P_{\sigma}=0.66666, P_{\mu}=1, k=-1$, (c) $\sigma_0 = 0, \mu_0=-0.001, D_{\sigma}=0.005, D_{\mu}=0.005, P_{\mu}=1, k=-1$, (d) $\sigma_0 = 0, \mu_0=-0.001, D_{\sigma}=0.005, D_{\mu}=0.005, P_{\sigma}=0.66666, k=-1$, and (e) $\sigma_0 = 0, \mu_0=-0.001, D_{\sigma}=0.005, D_{\mu}=0.005, P_{\sigma}=0.66666, P_{\mu}=1$.}
\label{fig:x_final plots}
\end{figure*}

Next, we separately plot the bifurcation diagrams for the last instances of $x_m$ against the parameters  $\sigma_0, \mu_0, P_{\sigma}, P_{\mu}$, and $k$. See Fig.~\ref{fig:x_final plots}. The first plot Fig.~\ref{fig:x_final plots}(a) depicts the bifurcation plot against $\sigma_0$. The parameters have been set to $P_{\mu} = 1, P_{\sigma} \sim 0.666, k=-1, \mu_0 = -0.001, D_{\mu} = 0.005$ and $D_{\sigma} = 0.005$. Keeping aside the values of $k$ and $\sigma_0$, note that we have utilized the same parameter combinations in Fig.~\ref{fig:sigma_0vskcolorplot}. 
Going back to the bifurcation diagram Fig.~\ref{fig:x_final plots}(a) we observe several cavities of synchronous states and asynchronous states almost evenly distributed in the range $-0.01<\sigma_0<0.01$. For example see the fourth plot in Fig.~\ref{fig:5plots1}(a) where we have used the same parameter combinations. Although there exists two clusters (one coherent, one solitary) and a single green node with $\Gamma_{2,m} \in [-0.15, 0.75]$, they are very close to each other, as evident from the $x-$scale in the $Y-$axis and obviously from the synchronization error value $E = 0.073$, which is quite low. 
Similar kind of patterns have been spotted by the authors in the references \cite{de2015effect}. identical arguments can be made for the bifurcation diagram where we vary $\mu_0$ in Fig.~\ref{fig:x_final plots}(b). Parameter values are $P_{\mu} = 1, P_{\sigma} \sim 0.666, k=-1, \sigma_0 = 0, D_{\mu} = 0.005$ and $D_{\sigma} = 0.005$.  
In the bifurcation plot shown in Fig.~\ref{fig:x_final plots}(b), we observe that for an approximate range of $\mu_0$ from $-0.00045$ to $-0.0002$, there is a gap where the nodes are synchronized. In the plot Fig.~\ref{fig:x_final plots}(e), we can thus state that the dynamics mostly exhibits asynchronity in the approximate range of $k$ from $0.01$ to around $3$. Outside this range, the dynamics portray synchronity. Similar inferences about Fig.~\ref{fig:x_final plots}(c), (d) can be drawn from the color-coded plots in Fig.~\ref{fig:sigma_0vsmu_0colorplot}.

\section{Conclusion}
\label{sec:conclusions}
In this article, we have made a substantial improvement in the ring-star network of Chialvo neurons under the influence of an electromagnetic flux by considering a heterogeneous topology of the network. The motivation behind this was to study the spatiotemporal dynamics of an ensemble of neurons which mimics a realistic nervous system. Heterogeneity is realized not only in space but also in time. We have introduced a noise modulation to incorporate heterogeneity in space and a time-varying structure of the links between neurons that update probabilistically. Noise sources are sampled from a uniform distribution. How the network dynamics change on induction of other complicated distributions opens a future scope of study. Exploring various computer-generated diagrams like phase portraits, spatiotemporal plots and recurrent plots, we observe rich dynamical behaviors like two-cluster states, solitary nodes, chimera states, traveling waves, and coherent states.\par
One of the purposes of the paper was to study the appearance of this special kind of asynchronous behavior called \textit{solitary nodes} and characterize them using a metric called \textit{cross-correlation coefficient}. It was observed that the cross-correlation coefficient indeed characterizes the solitary nodes in an efficient manner. Two more measures that we implemented here were the \textit{synchronization error} to study how synchronized the nodes behave, and \textit{sample entropy} to study the complexity of the network dynamics. Under different pairwise combinations of the model parameters we have studied the response of each of these three measures and have tried to conclude how efficiently they relate. Asynchronization and sample entropy, for example, portray a fairly convincing proportional relationship, whereas asynchronization and cross-correlation coefficient exhibit an inverse relation. These open a research direction to study how exactly they relate and if it is possible to establish an analytical relationship between them. Note that the trend obtained was a bit noisy. We hypothesize that this may be due to the heterogeneity introduced by noise modulations and time-varying links in the system. It would be interesting to see a noise-free relationship in a  noise-free system which we plan to address in the future.  Also, it would be very interesting to analytically and numerically study the existence of chaos in the ring-star network dynamics (by studying the Lyapunov exponent of the network system and methods to compute it).\par
Finally, we have also studied the one-parameter bifurcation diagrams of the synchronization of the nodes by plotting the final state of all the nodes against one of the parameters of the system. Interestingly, We observe some windows in the parameter regimes where the probability of all the nodes to be in a synchronized state is much higher than being asynchronous.\par
The fact that a ring-star network of chialvo neurons exhibits quite a rich dynamics, this study raises the question of how identical or contrasting the reported dynamical behaviors would be if we had considered different topologies and/or perturbations. Keeping this in mind, we plan to study the behaviors of an ensemble of Chialvo neurons under different topological arrangements: heterogeneous but quenched (no change in topology with time); has coupling strengths that are dependent on a normalized distance between a pair of neurons; is a multiplex heterogeneous network (See the recent work \cite{ShMu21b}); is perturbed by either temperature or photosensitivity. The fact that solitary nodes oscillate in a completely different phase from the main coherent ensemble, raises a future direction of exploring anti-phase synchronization with different coupling forms such as attractive and/or nonlinear repulsive couplings. A recent similar study \cite{ShMu21a} with Van der Pol oscillators has also been published.

\section{Acknowledgements}
The authors acknowledge Dr. David J. W. Simpson for his continuous guidance, support and sparing a part of his valuable time to go through the manuscript and suggest significant revisions. The authors also wish to thank Prof. Astero Provata, National Center for Scientific Research Demokritos and Aasifa Rounak, University College Dublin for their crucial inputs. Additionally, I.G acknowledges the mentoring provided by both S.S.M and H.O.F.

\bibliographystyle{plain}
\bibliography{main}

\end{document}